\documentclass[english]{article}
\usepackage[T1]{fontenc}
\usepackage[latin9]{inputenc}
\usepackage{geometry}
\usepackage{float}
\usepackage{amsmath}
\usepackage{amssymb}
\usepackage{graphicx}
\usepackage{xcolor}
\usepackage{algpseudocode}
\usepackage{algorithm}
\usepackage{scrextend}
\usepackage{pslatex}
\usepackage{graphicx}
\usepackage{amsthm}
\usepackage{amsfonts}
\usepackage{verbatim}
\usepackage{authblk}
\usepackage{lipsum}
\usepackage{background}
\usepackage[hidelinks]{hyperref}
\PassOptionsToPackage{normalem}{ulem}
\usepackage{ulem}
\makeatletter
\usepackage{placeins}
\usepackage{tikz}
\usetikzlibrary{decorations.pathreplacing,calc}
\newcommand{\tikzmark}[1]{\tikz[overlay,remember picture] \node (#1) {};}

\newcommand*{\AddNote}[4]{%
    \begin{tikzpicture}[overlay, remember picture]
        \draw [decoration={brace,amplitude=0.5em},decorate,ultra thick,red]
            ($(#3)!(#1.north)!($(#3)-(0,1)$)$) --  
            ($(#3)!(#2.south)!($(#3)-(0,1)$)$)
                node [align=center, text width=1.5cm, pos=0.5, anchor=west] {#4};
    \end{tikzpicture}
}%
\newcommand*{\AddNoteLeft}[4]{%
    \begin{tikzpicture}[overlay, remember picture]
        \draw [decoration={brace,amplitude=0.5em},decorate,ultra thick,red]
            ($(#3)!(#1.north)!($(#3)-(0,1)$)$) --  
            ($(#3)!(#2.south)!($(#3)-(0,1)$)$)
                node [align=center, text width=1.5cm, pos=0.5, anchor=east] {#4};
    \end{tikzpicture}
}%

\numberwithin{equation}{section}
\numberwithin{figure}{section}

\usepackage{xcolor}
\usepackage{algpseudocode}
\usepackage{algorithm}

\makeatletter

\theoremstyle{definition}

\theoremstyle{remark}

\numberwithin{equation}{section}
\providecommand{\tabularnewline}{\\}

\algnewcommand{\algorithmicgoto}{\textbf{go to}}%
\algnewcommand{\Goto}[1]{\algorithmicgoto~\ref{#1}}%

\makeatletter
\gdef\@ptsize{2}
\let\@currsize\normalsize
\makeatother
\usepackage{setspace}
\usepackage[font=small, width=0.9\textwidth]{caption}

\usepackage{microtype}  

\usepackage{footnote}
\makesavenoteenv{tabular}
\makesavenoteenv{table}

\usepackage{rotating}  
\usepackage{array}  
\usepackage{booktabs}  

\makeatletter

\providecommand{\tabularnewline}{\\}

\floatstyle{ruled}
\newfloat{algorithm}{tbp}{loa}
\providecommand{\algorithmname}{Algorithm}
\floatname{algorithm}{\protect\algorithmname}

\makeatother

\title{Stability and Complexity Analyses of Finite Difference
Algorithms for the Time-fractional Diffusion Equation}

\author[1]{Nirupama Bhattacharya}
\affil[1]{Department of Bioengineering\\University of California, San Diego}
\author[1,2,3]{Gabriel A. Silva \footnote{Corresponding author - G.S., Email: gsilva@ucsd.edu}}
\affil[2]{Department of Neurosciences\\University of California, San Diego}
\affil[3]{Center for Engineered Natural Intelligence\\University of California, San Diego}

\date{} 

\backgroundsetup{
  position=current page.north,
  angle=0,
  nodeanchor=north,
  vshift=-5mm,
  hshift=2mm,
  opacity=2,
  scale=3,
  contents=
}

\begin{document}


\maketitle
\renewenvironment{abstract}
{\begin{quote}
\noindent \rule{\linewidth}{.5pt}\par{\bfseries \abstractname.}}
{\medskip\noindent \rule{\linewidth}{.5pt}
\end{quote}
}

\begin{abstract}
Fractional differential equations (FDEs) are an extension of the theory of fractional calculus. However, due to the difficulty in finding analytical solutions, there have not been extensive applications of FDEs until recent decades. With the advent of increasing computational capacity along with advances in numerical methods, there has been increased interest in using FDEs to represent complex physical processes, where dynamics may not be as accurately captured with classical differential equations. The time-fractional diffusion equation is an FDE that represents the underlying physical mechanism of anomalous diffusion. But finding tractable analytical solutions to FDEs is often much more involved than solving for the solutions of integer order differential equations, and in many cases it is not possible to frame solutions in a closed form expression that can be easily simulated or visually represented. Therefore the development of numerical methods is vital. In previous work we implemented the full 2D time-fractional diffusion equation as a Forward Time Central Space finite difference equation by using the Gr\"{u}nwald-Letnikov definition of the fractional derivative. In addition, we derived an adaptive time step version that improves on calculation speed, with some tradeoff in accuracy. Here, we explore and characterize stability of these algorithms, in order to define bounds on computational parameters that are crucial for performing accurate simulations. We also analyze the time complexity of the algorithms, and describe an alternate adaptive time step approach that utilizes a linked list implementation, which yields better algorithmic efficiency.   \\
\end{abstract} 



\section{Introduction}

Fractional calculus as a mathematical field has existed almost as
long as classical integer-order calculus, much of the theory having
been developed by Liouville, Laplace, Fourier, Euler, Lagrange,
Riemann, and others \cite{Oldham_Spanier, Magin:2006uo}. Analogous to ordinary and partial differential
equations being used to model physical phenomena using classical calculus,
fractional differential equations (FDEs) are an extension of the theory
of fractional calculus. However, due to the difficulty in finding
analytical solutions to FDEs, there have not been extensive applications
of these equations until recent decades. With the advent of increasing
computational capacity along with advances in numerical methods, there has been increased interest in using FDEs to represent complex physical processes, where dynamics may not be as accurately captured with
classical differential equations.

Fractional differential equations have been useful for modeling such
varied applications as advection-dispersion transport
phenomena \cite{Baeumer_adv_dispersion}, wave propagation in bone and other rigid porous
materials \cite{Sebaa_Fellah2006, Fellah_2002}, viscoelasticity \cite{Soczkiewicz_Viscoelasticity}, 
edge detection in image processing \cite{Mathieu2003}, and neurodynamics and modeling of signal
processing in neuronal dendrites \cite{Henry_Langlands_2008a}. Richard Magin discusses
a wide array of fractional calculus modeling applications in the field
of bioengineering alone, including signaling across membranes,
feedback control in neural systems, neurodynamics, capacitor and dielectric
models, the behavior of viscoelastic materials, circuit models of
electrode interfaces, cell biomechanics, and electrochemistry, among many other applications \cite{Magin:2006uo, Magin_Ovadia_2008}.

In particular, the time-fractional diffusion equation is an FDE that represents the 
underlying physical mechanism of anomalous diffusion. In the same
way that the classical diffusion equation can be derived from statistical
analysis of particle interactions, one can also show how the micromolecular
behavior of particles can, under a different set of fundamental statistical
assumptions, be represented on a macro scale in two dimensions by the fractional diffusion
equation of the form \cite{Thesis}: 
\begin{equation}
\frac{\partial^{\gamma}u(\vec{x},t)}{\partial t^{\gamma}}=\alpha\frac{\partial^{2}u(\vec{x},t)}{\partial x^{2}}+\beta\frac{\partial^{2}u(\vec{x},t)}{\partial y^{2}}\label{eq:frac-diffn-eqn}
\end{equation}
where $\gamma<1$ denotes the subdiffusion regime, $\gamma=1$ denotes
classical Gaussian diffusion, $\gamma>1$ denotes the superdiffusion
regime, and $\alpha$ and $\beta$ represent the diffusion coefficients
in the $x$ and $y$ directions, respectively (in $\frac{spatial\; unit^2}{time\; unit^\gamma}$).

Finding tractable analytical solutions to FDEs like Eq. \ref{eq:frac-diffn-eqn}
is often much more involved than solving for the solutions of integer
order differential equations, and in many cases it is not possible
to frame solutions in a closed form expression that can be easily
simulated or visually represented. Therefore the development of numerical
methods is vital. Over the last decade there have been numerous 
algorithms developed for FDEs like the time-fractional diffusion equation,
and there is often a tradeoff between calculation speed and efficiency,
complexity, stability, and accuracy. In \cite{Paper0} we walk
through the discretization of Eq. \ref{eq:frac-diffn-eqn} into the
full two-dimensional fractional FTCS (Forward Time Central Space) finite difference equation, making
use of the Gr\"{u}nwald-Letnikov definition of the fractional derivative,
which allows us to use a discrete summation in our numerical algorithm.
We also derive an adaptive time step algorithm
which builds on the foundation of the full 2D fractional FTCS equation
but improves calculation speed and efficiency, while maintaining accuracy \cite{Paper0}. However, both the FTCS equation and adaptive algorithms
are explicit methods with a limited stability regime. Here we fully
explore and characterize the stability of these two finite difference
schemes in order to define bounds on important computational parameters like
time step and spatial discretization and grid size; selecting these
parameters appropriately is crucial for performing accurate simulations.

There are numerous approaches to analyzing the stability of finite
difference schemes, including modified wavenumber analysis, matrix
eigenvalue analysis, and other more mathematically complex methods
derived from stability definitions involving matrix norms. Several
methods, including matrix analysis, are widely applicable but sometimes
impossible to approach analytically or without the aid of iterative
procedures. However, since our numerical algorithms are
both fully discretized in space and time, and linear with constant
coefficients (in their homogeneous versions), we show in the
next section that an appropriate choice is a von Neumann stability analysis.
We assume the general solution at an arbitrary location and time step
is a linear combination of special solutions, represented
by a function that is separable in its temporal and spatial variables.
We find that our stability analysis and expressions of parameter bounds
agree very well with our results, where we simulate several examples
with various parameter combinations.

In addition to stability analysis, a complexity analysis is another
important metric used to characterize the efficiency of numerical
algorithms, as measured by execution time of the algorithm as a function
of some input variable. In the case of the algorithms described in
this paper, we are interested in how the execution time varies
with $N_{x},\; N_{y}$(grid size of the simulation in the x and y
directions, respectfully), or $N$, the number of timesteps. The relationship
between execution time and input variable is usually given in Big-$O$
notation using a worst case scenario, or average
case scenario which is often a better reflection of the average behavior
of the algorithm run time. In Section 3, we explore the time complexity of the same two algorithms
for which we complete the stability analysis, and in addition, analyze an alternate
version of the adaptive timestep algorithm that involves a linked
list implementation (introduced in \cite{Paper0}) that yields better algorithmic efficiency. For
mathematical simplicity we analyze our
algorithms according to worst case scenarios and interpret the results as a proof of bounding behavior
- the actual run time of the algorithms may be more efficient, but
never less efficient. We also present simulated data that verifies our theoretical
complexity analyses.

\section{Stability Analysis of the Two-dimensional Fractional FTCS Discretization}

\subsection{Full Implementation}

We begin by considering the circumstances under which the 2D fractional
FTCS finite difference equation is stable, and consider
the advantages and disadvantages of several different analytical approaches.
We restate the FTCS discretized equation in two dimensions (see \cite{Paper0, Thesis} for more details): 
\begin{eqnarray}
u_{j,l}^{n+1}-u_{j,l}^{n} & = & \Delta_{t}^{\gamma}\sum_{m=0}^{n}\psi(\gamma,m)\left(\frac{\alpha}{\Delta x^{2}}\delta x_{j,l}^{n-m}+\frac{\beta}{\Delta y^{2}}\delta y_{j,l}^{n-m}\right)\label{eq:2d discretized full equation}\\
\delta x_{j,l}^{n} & = & u_{j+1,l}^{n}-2u_{j,l}^{n}+u_{j-1,l}^{n} \nonumber \\
\delta y_{j,l}^{n} & = & u_{j,l+1}^{n}-2u_{j,l}^{n}+u_{j,l-1}^{n} \nonumber
\end{eqnarray}
where $\Delta_{t}$ is the time step, and $\Delta x$ and $\Delta y$
are spatial grid discretizations in the x and y directions. Our spatial range is

\[
u_{j,l}^{n}\mid j=0,1,2,...,N_{x}-1;\; l=0,1,2,...,N_{y}-1
\]
and 
\[
u_{j,l}^{n}\mid j=0,N_{x}-1 ;\; l=0,N_{y}-1
\]
denotes the boundary points. $\psi\left(\gamma,m\right)$, which we
refer to as a `memory function', represents a binomial term $\left(-1\right)^{m}\left(\begin{array}{c}
1-\gamma\\
m
\end{array}\right)$, and results from the use of the Gr\"{u}nwald-Letnikov definition of
the fractional derivative, as discussed in \cite{Paper0,podlubny1999fractional,Thesis}.

\subsubsection{Matrix Stability Analysis}

Matrix stability analysis is a method of analysis that can be applied
to any problem without particular assumptions or constraints on coefficients,
and includes the effects of boundary conditions \cite{Parviz}. We would like to
reformulate Eq. \ref{eq:2d discretized full equation} into a matrix
equation in the form $W_{n+1}=L_{n}W_{n}+f_{n}$ where the $W$ matrices
hold $u$ values at times $n+1$ and $n$, the $L_{n}$ matrix represents
the difference operator at time $n$, and $f_{n}$ is a constant vector that takes
into account boundary conditions at time step $n$ (here we assume
boundary conditions are known for all timesteps).

For example, for $n=0$, Eq. \ref{eq:2d discretized full equation} becomes: 
\[
u_{j,l}^{1}=u_{j,l}^{0}+\Delta_{t}^{\gamma}\psi(\gamma,0)\left(\frac{\alpha}{\Delta x^{2}}\delta x_{j,l}^{0}+\frac{\beta}{\Delta y^{2}}\delta y_{j,l}^{0}\right)
\]
Then the matrix equation is formulated as follows. We order grid values $u$ into a vector (excluding the boundary points) such that $W$ at
some time step $n$ is 
\begin{equation}
W_{n}=\left[\begin{array}{c}
u_{1,1}\\
u_{1,2}\\
u_{1,3}\\
\vdots\\
u_{1,N_{y}-2}\\
\vdots\\
u_{N_{x}-2,1}\\
u_{N_{x}-2,2}\\
\vdots\\
u_{N_{x}-2,N_{y}-2}
\end{array}\right]^{n}\label{eq:grid_layout}
\end{equation}

\noindent For $n = 0$ our matrix equation now becomes $W_{1}=L_{0}W_{0}+f_{0}$
where the matrix $L_{0}$ accounts for interior non-boundary grid points and
consists of block matrices $T_{0}$ (a tridiagonal matrix) and $R_{0}$ (a
diagonal matrix): 
{\footnotesize
\begin{eqnarray*}
L_{0} & = & \left[\begin{array}{ccccc}
		T_{0} & R_{0} & 0 & \cdots & 0\\
		R_{0} & T_{0} & \ddots & \ddots & \vdots\\
		0 & R_{0} & \ddots & R_{0} & 0\\
		\vdots & \ddots & \ddots & T_{0} & R_{0}\\
		0 & \cdots & 0 & R_{0} & T_{0}
		\end{array}\right]\\
T_{0} & = & \left[\begin{array}{ccccc}
		1-2r_{x}\psi\left(\gamma,0\right)-2r_{y}\psi\left(\gamma,0\right) & r_{y}\psi(\gamma,0) & 0 & \cdots & 0\\
		r_{y}\psi(\gamma,0) &  & \ddots & \ddots & \vdots\\
		0 & \ddots & \ddots &  & 0\\
		\vdots & \ddots &  &  & r_{y}\psi(\gamma,0)\\
		0 & \cdots & 0 & r_{y}\psi(\gamma,0) & 1-2r_{x}\psi\left(\gamma,0\right)-2r_{y}\psi\left(\gamma,0\right)
		\end{array}\right]\\
R_{0} & = & \left[\begin{array}{cccc}
		r_{x}\psi(\gamma,0) & 0 & \cdots & 0\\
		0 & \ddots & \ddots & \vdots\\
		\vdots & \ddots & \ddots & 0\\
		0 & \cdots & 0 & r_{x}\psi(\gamma,0)
		\end{array}\right]
\end{eqnarray*}}
where $r_{x}=\frac{\alpha\Delta_{t}^{\gamma}}{\Delta x^{2}}$ and
$r_{y}=\frac{\beta\Delta_{t}^{\gamma}}{\Delta y^{2}}$. $L_{0}$ is size $(N_{x}-2) (N_{y}-2)$ by $(N_{x}-2)(N_{y}-2)$. $T_{0}$ and $R_{0}$ are each of size $(N_{y}-2)$ by $(N_{y}-2)$. $f_{0}$ takes into account the boundary points and is a column vector of constants.

\noindent For $n=1$ we now must take into account multiple past time
points, and our matrix equation becomes 
{\footnotesize
\begin{eqnarray*}
W_{2} & = & L_{1}W_{0}+L_{0}W_{1}+f_{01}\\
L_{1} & = & \left[\begin{array}{ccccc}
		T_{1} & R_{1} & 0 & \cdots & 0\\
		R_{1} & T_{1} & \ddots & \ddots & \vdots\\
		0 & R_{1} & \ddots & R_{1} & 0\\
		\vdots & \ddots & \ddots & T_{1} & R_{1}\\
		0 & \cdots & 0 & R_{1} & T_{1}
		\end{array}\right]\\
T_{1} & = & \left[\begin{array}{ccccc}
		-2r_{x}\psi\left(\gamma,1\right)-2r_{y}\psi\left(\gamma,1\right) & r_{y}\psi(\gamma,1) & 0 & \cdots & 0\\
		r_{y}\psi(\gamma,1) &  & \ddots & \ddots & \vdots\\
		0 & \ddots & \ddots &  & 0\\
		\vdots & \ddots &  &  & r_{y}\psi(\gamma,1)\\
		0 & \cdots & 0 & r_{y}\psi(\gamma,1) & -2r_{x}\psi\left(\gamma,1\right)-2r_{y}\psi\left(\gamma,1\right)
		\end{array}\right]\\
R_{1} & = & \left[\begin{array}{cccc}
		r_{x}\psi(\gamma,1) & 0 & \cdots & 0\\
		0 & \ddots & \ddots & \vdots\\
		\vdots & \ddots & \ddots & 0\\
		0 & \cdots & 0 & r_{x}\psi(\gamma,1)
		\end{array}\right]
\end{eqnarray*}
}
\normalsize

and where $f_{01}$ is a column vector that takes into account the boundary conditions from timesteps
$0$ and $1$. We can repurpose this multistep matrix equation into
a single step equation: 
\begin{eqnarray*}
W_{combined,2} & = & L_{combined,1}W_{combined,1} + f_{01}\;\mbox{with}\\
W_{combined,2} & = & \left[\begin{array}{c}
		W_{2}\\
		W_{1}
		\end{array}\right],\; W_{combined,1}=\left[\begin{array}{c}
		W_{1}\\
		W_{0}
		\end{array}\right],\; L_{combined,1}=\left[\begin{array}{cc}
		L_{0} & L_{1}\\
		I & 0
		\end{array}\right]
\end{eqnarray*}

and $I$ is simply the identity matrix. For $n=2$ the complexity
of the matrices involved continues to grow with 

\begin{eqnarray*}
W_{combined,3} & = & L_{combined,2}W_{combined,2}+f_{012}\\
\left[\begin{array}{c}
W_{3}\\
W_{2}\\
W_{1}
\end{array}\right] & = & \left[\begin{array}{ccc}
		L_{0} & L_{1} & L_{2}\\
		I & 0 & 0\\
		0 & I & 0
		\end{array}\right]\left[\begin{array}{c}
		W_{2}\\
		W_{1}\\
		W_{0}
		\end{array}\right]+f_{012}
\end{eqnarray*}
where $L_{0},L_{1},L_{2}...$ are block matrices which consist of
tridiagonal and diagonal sub-matrices, as demonstrated in previous
steps. In general, we have 

\begin{align}
W_{combined,n+1}&= L_{combined,n}W_{combined,n}+f_{0...n} \label{eq:gen-matrix-formula} \\
W_{combined,n+1} & = \left[\begin{array}{c}
		W_{n+1}\\
		\vdots\\
		\vdots\\
		\vdots\\
		W_{1}
		\end{array}\right] \nonumber \\
L_{combined, n} & = 
		\left[\begin{array}{ccccc}
		L_{0} & L_{1} & \cdots & L_{n-1} & L_{n}\\
		I & 0 &  & \cdots & 0\\
		0 & \ddots & \ddots &  & \vdots\\
		\vdots & \ddots & \ddots & \ddots\\
		0 & \cdots & 0 & I & 0
		\end{array}\right] \nonumber \\
L_{n} & = \left[\begin{array}{ccccc}
		T_{n} & R_{n} & 0 & \cdots & 0\\
		R_{n} & T_{n} & \ddots & \ddots & \vdots\\
		0 & R_{n} & \ddots & R_{n} & 0\\
		\vdots & \ddots & \ddots & T_{n} & R_{n}\\
		0 & \cdots & 0 & R_{n} & T_{n}
		\end{array}\right]\nonumber \\
T_{n} & = \left[\begin{array}{ccccc}
		-2r_{x}\psi\left(\gamma,n\right)-2r_{y}\psi\left(\gamma,n\right) & r_{y}\psi(\gamma,n) & 0 & \cdots & 0\\
		r_{y}\psi(\gamma,n) &  & \ddots & \ddots & \vdots\\
		0 & \ddots & \ddots &  & 0\\
		\vdots & \ddots &  &  & r_{y}\psi(\gamma,n) \nonumber \\
		0 & \cdots & 0 & r_{y}\psi(\gamma,n) & -2r_{x}\psi\left(\gamma,n\right)-2r_{y}\psi\left(\gamma,n\right)
		\end{array}\right], \text{for } n \geq 1\\
R_{n} & = \left[\begin{array}{cccc}
		r_{x}\psi(\gamma,n) & 0 & \cdots & 0\\
		0 & \ddots & \ddots & \vdots\\
		\vdots & \ddots & \ddots & 0 \nonumber \\
		0 & \cdots & 0 & r_{x}\psi(\gamma,n)
		\end{array}\right]
\end{align}

The goal with matrix stability analysis is to show that for a matrix
problem Eq. \ref{eq:gen-matrix-formula}, the spectral radius $\rho(L_{combined,n})$
of the difference matrix $L_{combined,n}$, (i.e., the modulus of its largest eigenvalue)
is bounded, usually by 1 \cite{Polytech_course_notes}. This involves finding the maximum eigenvalue
for $L_{combined,n}$, for all timesteps $n$, and proving that its modulus is
$\leq1$. The matrix $L_{combined,n}$ is quite complex
and is not in the form of a special matrix where there are known analytical
algorithms to find the eigenvalues. This makes the process of finding
the eigenvalues very complicated and most likely beyond the realm of an analytic solution, so we consider this
to be outside the scope of this paper, and conclude that despite the
difficulty of finding the spectral radius, a major advantage of this
method would be its generality - it could be applied to any set of boundary conditions,
or schemes with non-constant coefficients.

\subsubsection{Von Neumann Analysis}
\label{subsection:vonNeumann}

As we saw in the previous section, matrix stability analysis for complex
algorithms with no simplifying assumptions, requires methods for finding
eigenvalues, possibly including iterative numerical schemes. However, our finite difference scheme in Eq. \ref{eq:2d discretized full equation}
has the advantage of being linear with constant coefficients (i.e.
$\alpha,\beta$ don't depend on $t$ or $\vec{x}$) on a uniformly spaced grid. With an additional
simplifying assumption of periodic boundary conditions, we can make
use of von Neumann stability analysis and Fourier series expansion.
While von Neumann analysis doesn't take into account boundary conditions
in the same way matrix stability analysis does, we make the reasonable
assumption that stability issues are mostly affected by discretization of differential
equations inside the domain, and minimally affected by boundary conditions \cite{Parviz}.
In the majority of modeling and simulation problems in engineering,
we consider this to be a valid and applicable simplifying assumption.

Given our linear scheme with constant coefficient and periodic boundary
conditions, we assume the solution of scheme \ref{eq:2d discretized full equation}
at a particular location and time point, can be expressed as a discrete
finite complex Fourier series expansion: 
\begin{eqnarray}
u_{j,l}^{n} & = & \sum_{d=0}^{N_{y}-1}\sum_{c=0}^{N_{x}-1}A_{cd}^{n}e^{ik_{c}x_{j}}e^{ik_{d}y_{l}}\label{eq:von-Neumann-F-expansion}
\end{eqnarray}
where  $A_{cd}^n$ are the Fourier mode amplitudes, $x_{j}=j\Delta x,\; j=0,1,...,N_{x}-1, y_{l}=l\Delta y,\; l=0,1,...,N_{y}-1$, and 
$k_{c}$ and $k_{d}$ are the spatial wavenumbers in the $x$ and
$y$ directions, respectively. The wavenumbers represent spatial frequency
and are related to wavelength by $\lambda_{x}=\frac{2\pi}{k_{c}},\;\lambda_{y}=\frac{2\pi}{k_{d}}$.
In the context of Fourier modes and stability analysis, we are interested
in suppressing (or as we shall see later in this section, limiting
the amplification factor of) the high frequency modes (represented
by wavelengths on the order of grid sizes $\Delta x$ and $\Delta y$,
and wavenumber values away from $0$) because these modes
have a low probability of corresponding to any real features of the
true solution, and usually represent noise arising from finite precision
arithmetic for numerical calculations. The time evolution of each Fourier mode is determined by the same time marching scheme as the full solution $u_{j,l}^{n}$. And since no Fourier mode can be allowed to increase unbounded in time, we can single out any arbitrary harmonic $A_{cd}^{n}e^{ik_{c}x_{j}}e^{ik_{d}y_{l}}$ and substitute this form into Eq. \ref{eq:2d discretized full equation} \cite{num_methods_flows}. We utilize the relationships $x_{j+1}=x_{j}+\Delta x,\; x_{j-1}=x_{j}-\Delta x,\; y_{l+1}=y_{l}+\Delta y\;\mbox{and \ensuremath{y_{l-1}=y_{l}-\Delta y}}$, as well as Euler's formula and half angle trigonometric identities,
to obtain 

\begin{align}
A_{cd}^{n+1}-A_{cd}^{n}+\Delta_{t}^{\gamma}\left(\frac{4\alpha}{\Delta x^{2}}sin^{2}\left(\frac{k_{c}\Delta x}{2}\right)+\frac{4\beta}{\Delta y^{2}}sin^{2}\left(\frac{k_{d}\Delta y}{2}\right)\right)\sum_{m=0}^{n}\psi\left(\gamma,m\right)A_{cd}^{n-m}=0
\end{align}

	If we make the assumption that $A_{cd}$ is independent of time step,
we can write 
\begin{equation}
A_{cd}^{n+1}=\sigma_{cd}A_{cd}^{n}\label{eq:amp_factor_eqn}
\end{equation}
where $\sigma_{cd}$ is known as the amplification factor that gives
the growth or decay for Fourier mode $\left(c,d\right)$. We require
$\left|\sigma_{cd}\right|\leq1$ for all Fourier modes, or else some
modes will be amplified at every time step and dominate the solution.
We divide Eq. \ref{eq:amp_factor_eqn} by $A_{cd}^{n}$ and rearrange to get 
\[
\sigma_{cd}=1-\Delta_{t}^{\gamma}\left(\frac{4\alpha}{\Delta x^{2}}sin^{2}\left(\frac{k_{c}\Delta x}{2}\right)+\frac{4\beta}{\Delta y^{2}}sin^{2}\left(\frac{k_{d}\Delta y}{2}\right)\right)\sum_{m=0}^{n}\psi\left(\gamma,m\right)\sigma_{cd}^{-m}
\]
We require $\left|\sigma_{cd}\right|\leq1$ for all values of $\sigma_{cd}$,
which are the roots of this polynomial. We can also see that the difficulty
of solving for the roots increases as time progresses and $n$ increases.
Past the first few timesteps, solving for the roots of the polynomial requires
iterative matrix procedures and does not result in a tractable analytical
solution from which we can find the bounds on the relationship between
$\alpha,\beta,\Delta t,\Delta x,\mbox{and }\Delta y$. Instead of solving
for roots, we will approach the problem by considering `worst-case'
scenarios to simplify our expression, and proceed in a similar manner
as Yuste and Acedo's treatment of the analogous one-dimensional case \cite{yuste2003explicit}. Consider the bound on the parameters resulting from the
case $\sigma_{cd}=1$: 
\begin{equation}
0\geq\Delta_{t}^{\gamma}\left(\frac{4\alpha}{\Delta x^{2}}sin^{2}\left(\frac{k_{c}\Delta x}{2}\right)+\frac{4\beta}{\Delta y^{2}}sin^{2}\left(\frac{k_{d}\Delta y}{2}\right)\right)\sum_{m=0}^{n}\psi\left(\gamma,m\right)\label{eq:sigma=1}
\end{equation}
The summation term and parameters are always greater than $0$, and
therefore Eq. \ref{eq:sigma=1} can only hold true if both
wavenumbers are $0$ (a constant function). Therefore, we consider the other `worst-case'
scenario, where $\sigma_{cd}=-1$: 
\begin{equation}
2=\Delta_{t}^{\gamma}\left(\frac{4\alpha}{\Delta x^{2}}sin^{2}\left(\frac{k_{c}\Delta x}{2}\right)+\frac{4\beta}{\Delta y^{2}}sin^{2}\left(\frac{k_{d}\Delta y}{2}\right)\right)\sum_{m=0}^{n}\psi(\gamma,m)\left(-1\right)^{m}\label{eq:sigma=00003D00003D-1}
\end{equation}
From Eq. \ref{eq:sigma=00003D00003D-1} we can infer the inequality
bounding the parameters as 
\begin{equation}
\Delta_{t}^{\gamma}\left(\frac{4\alpha}{\Delta x^{2}}sin^{2}\left(\frac{k_{c}\Delta x}{2}\right)+\frac{4\beta}{\Delta y^{2}}sin^{2}\left(\frac{k_{d}\Delta y}{2}\right)\right)\leq\frac{2}{\sum_{m=0}^{n}\psi(\gamma,m)\left(-1\right)^{m}}=B_{1}\left(\gamma,n\right)\label{eq:B1-full}
\end{equation}
We can see that the bounded value is dependent on the time step $n$
but the dependence is very weak. Yuste and Acedo demonstrate that in the limit $n\rightarrow\infty$, the summation term in
Eq. \ref{eq:B1-full} involving the memory function, is an alternating
converging series; as $n$ increases, $B_{1}$ oscillates but clearly
converges to an equilibrium value \cite{yuste2003explicit}. Therefore
we can approximate $ $$\sum_{m=0}^{n}\psi(\gamma,m)\left(-1\right)^{m}$
by taking the limit as $n\rightarrow\infty$. In \cite{Paper0} we originally derived
the 2D FTCS full finite difference equation and adaptive memory algorithms
by using the Gr\"{u}nwald-Letnikov definition to make a first-order approximation
of the Riemann-Liouville fractional derivative operator \cite{podlubny1999fractional}.
And in the first-order approximation, the memory function $\psi(\gamma,m)$
is defined as $\left(-1\right)^{m}\left(\begin{array}{c}
1-\gamma\\
m
\end{array}\right)$. Podlubny (\cite{podlubny1999fractional}) shows that this expression can also be defined as the coefficients
of the power series for the function $\left(1-z\right){}^{1-\gamma}:$
\begin{equation}
\left(1-z\right)^{1-\gamma}=\sum_{m=0}^{\infty}\left(-1\right)^{m}\left(\begin{array}{c}
1-\gamma\\
m
\end{array}\right)z^{m}=\sum_{m=0}^{\infty}\psi(\gamma,m)z^{m}\label{eq:power series coefficients-1}
\end{equation}
Applying this identity to $B_{1}\left(\gamma,n\right)$ and considering
the limit $n\rightarrow\infty$, we can derive 
\begin{equation}
B_{2}\left(\gamma\right)=2^{\gamma}\label{eq:B2-full}
\end{equation}

As Yuste and Acedo describe in their considerations of stability in
the one-dimensional case, we could consider the second-order approximation of the
Riemann-Liouville definition by making use of a different set coefficients
to define our memory function. This option would result in a slightly
lower bound on our parameters, but for the purposes of estimating
appropriate parameter values, we consider the first-order approximation
to be sufficient.

From (\ref{eq:B1-full}) we can see that we have three degrees of
freedom in parameters that we are free to choose for our simulations
: $\Delta_{t},\;\Delta x^{2},\;\Delta y^{2}$ . To determine the most
conservative bounds on those degrees of freedom (minimum values
for grid size and maximum value for time step), we assume the maximum
possible value for the trigonometric terms, $sin^{2}\left(\frac{k_{c}\Delta x}{2}\right)=sin^{2}\left(\frac{k_{d}\Delta y}{2}\right)=1$.
The resulting inequality determining parameter bounds is 
\[
\Delta_{t}^{\gamma}\left(\frac{4\alpha}{\Delta x^{2}}+\frac{4\beta}{\Delta y^{2}}\right)\leq B_{2}\left(\gamma\right)
\]
If we assume a simple case with a square grid ($\Delta x=\Delta y=\Delta)$
and equal diffusion coefficients in both spatial directions $\left(\alpha=\beta\right)$,
the expression becomes

\noindent 
\begin{equation}
r=\frac{\alpha\Delta t^{\gamma}}{\Delta^{2}}\leq\frac{B_{2}(\gamma)}{8}=2^{\gamma-3}=B_{full}\left(\gamma\right)\label{eq:Bound-full}
\end{equation}

As a final check, when we set $\gamma=1$ (classical diffusion),
we recover the well known traditional bound $r\leq\frac{1}{4}$ that
results from the classical 2D diffusion equation using the FTCS discretization.

\subsubsection{Error Propagation}

An alternative approach to stability is to consider roundoff error from floating point or finite precision arithmetic, and ensure that this does not propagate in time by requiring error $\left|\epsilon_{j,l}^{n+1}\right| \leq \left|\epsilon_{j,l}^{n}\right|$. A mathematical derivation based on this interpretation would progress in much the same manner as in Section \ref{subsection:vonNeumann}
and result in similar parameter boundaries.

\subsection{Adaptive Memory Algorithm}

As discussed in \cite{Paper0}, while the full implementation
represented by Eq. \ref{eq:2d discretized full equation} is an accurate
finite difference approximation to the solution of Eq. \ref{eq:frac-diffn-eqn},
using the Gr\"{u}nwald-Letnikov definition of the fractional derivative
involves a summation that takes into account the entire past history
of the simulation. As the simulation progresses, calculating the summation
term becomes increasingly cumbersome, time consuming, and memory intensive.
Therefore we have developed an adaptive memory method (introduced in \cite{Paper0}) that improves on computational efficiency while maintaining
a level of accuracy comparable to the original full implementation.
This is achieved by recognizing that the further back a time point
is in the history of the simulation, the less it contributes to the calculation
of the solution at the next time point. Therefore we sample the history
of the system more frequently for recent time points (which contribute
more to the solution at the next time step) and less often for time
points further back in the history of the system, weighted by an appropriate
amount to compensate for less frequent sampling; this approach significantly
reduces actual computational time. For smooth functions, the convolution
of $\psi\left(\gamma,m\right)$ and the function in question changes slowly enough that the weighting of the sampled time points compensates
for the less frequent sampling, and maintains overall accuracy. The
motivation and derivation of this adaptive algorithm is discussed in further
detail in \cite{Paper0}. Here we reproduce a modified version of the algorithm and recognize that
analyzing stability requires only a few modifications from the process
described in the last section.

\noindent 
\begin{eqnarray}
\frac{u_{j,l}^{n+1}-u_{j,l}^{n}}{\Delta_{t}^{\gamma}} & = & \sum_{m=0}^{a}
\overbrace{\psi(\gamma,m)\left(\frac{\alpha}{\Delta x^{2}}\delta x_{j,l}^{n-m}+\frac{\beta}{\Delta y^{2}}\delta y_{j,l}^{n-m}\right)}^{R1}\label{eq:adaptive-memory}\\
 &  & +\sum_{s=2}^{s_{max}(n)}\left\{ \sum_{\eta=1}^{\eta_{max}\left(s,n\right)} \overbrace{\left(2s-1\right)\psi\left(\gamma,M\left(s,\eta\right)\right)\left(\frac{\alpha}{\Delta x^{2}}\delta x_{j,l}^{n-M(s,\eta)}+\frac{\beta}{\Delta y^{2}}\delta y_{j,l}^{n-M(s,\eta)}\right)}^{R2}\right.\nonumber \\
 &  & +\left.\sum_{p=M(s,\eta_{max}) +s}^{min(a^{s},n)}
 \overbrace{\psi(\gamma,p)\left(\frac{\alpha}{\Delta x^{2}}\delta x_{j,l}^{n-p}+\frac{\beta}{\Delta y^{2}}\delta y_{j,l}^{n-p}\right)}^{R3}\right\} \nonumber \\
M\left(s,\eta\right) & = & a^{s-1}+\left(2s-1\right)\eta-s+1\nonumber \\
\eta_{max}\left(s,n\right) & = & min\left(\left\lfloor \frac{a^{s}-a^{s-1}}{2s-1}\right\rfloor ,\left\lfloor \frac{n-a^{s-1}}{2s-1}\right\rfloor \right)\nonumber \\
\delta x_{j,l}^{n} & = & u_{j+1,l}^{n}-2u_{j,l}^{n}+u_{j-1,l}^{n}\nonumber \\
\delta y_{j,l}^{n} & = & u_{j,l+1}^{n}-2u_{j,l}^{n}+u_{j,l-1}^{n}
\end{eqnarray}
where $a$ is the predefined base interval, $\left\lfloor \;\right\rfloor $
denotes the floor function, and $s_{max}$ is determined by the current
time step $n$ such that $a^{s_{max}-1}+1\leq n\leq a^{s_{max}}$. The terms $R1, R2, R3$ are referred to in more detail during the complexity analysis in later sections.

\begin{figure}
\begin{center}
a)\includegraphics[width=4.8in]{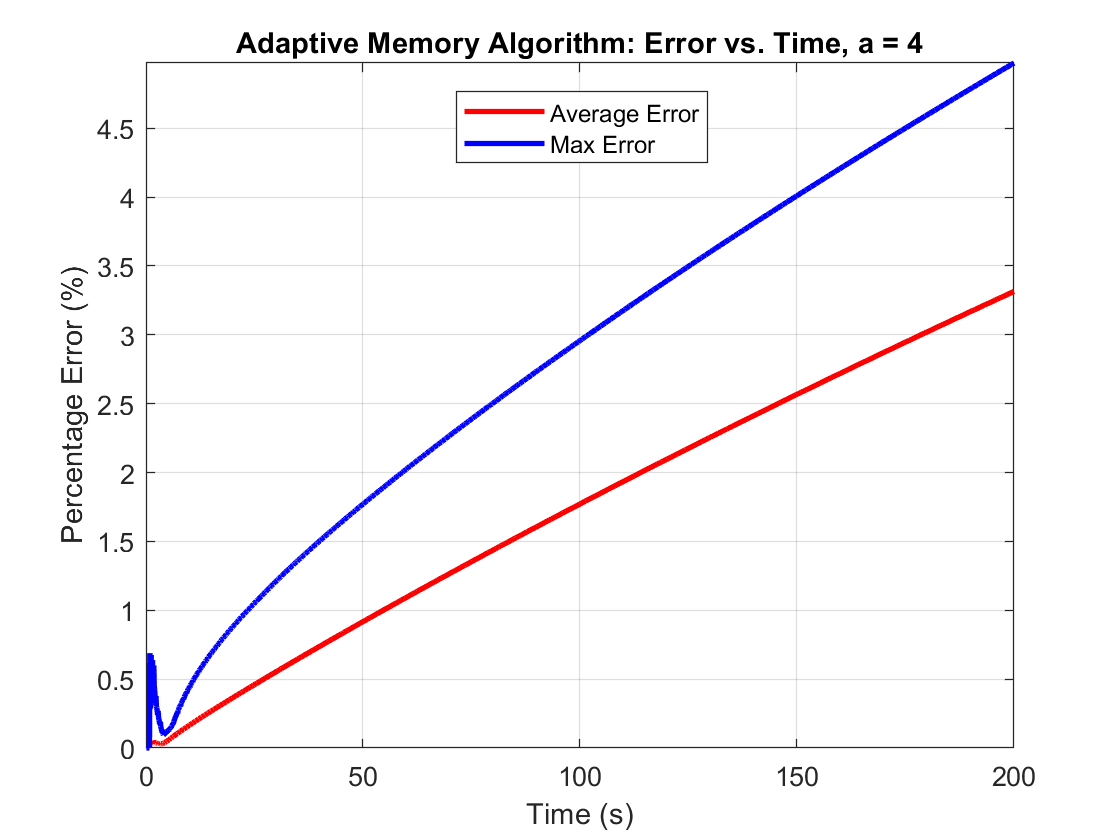}\\
b)\includegraphics[width=4.8in]{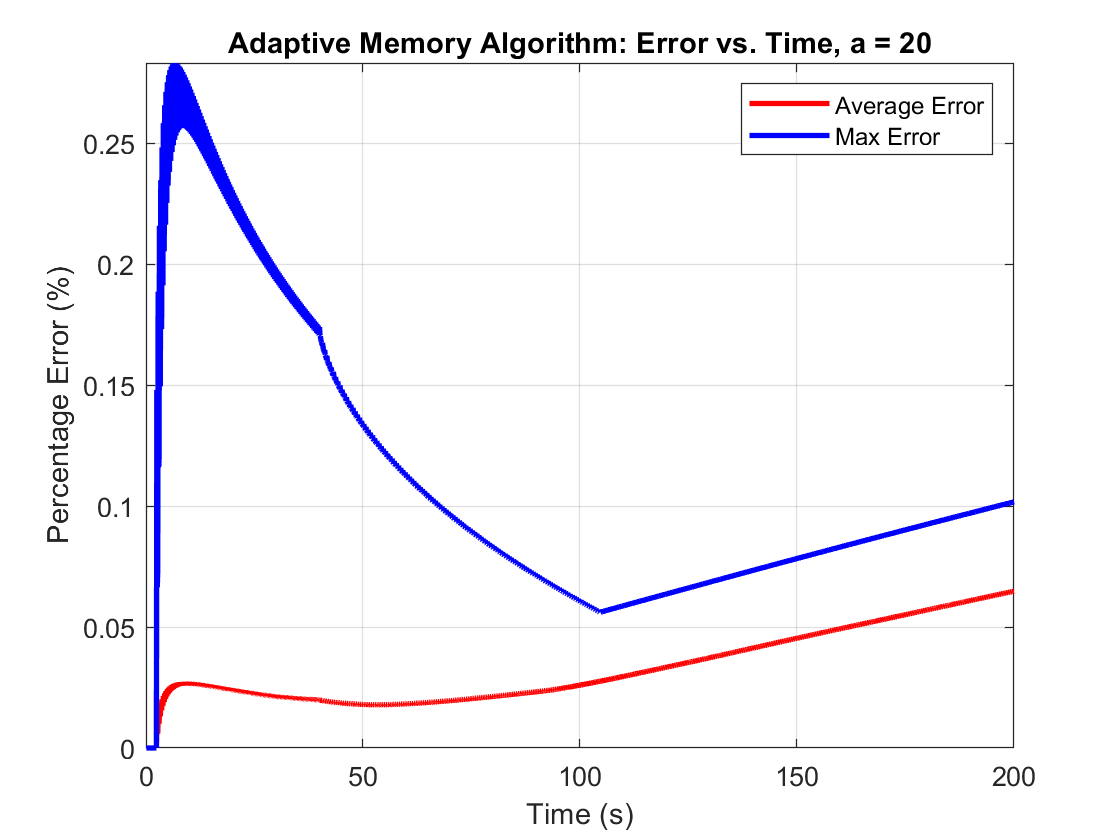} 
\caption[Effect of $a$ on Accuracy]{\footnotesize
\textbf{Effect of $a$ on Accuracy.} \textbf{a)} With adaptive step
parameter $a=4$, we see the error
(compared against the full 2D discretization defined in Eq. \ref{eq:2d discretized full equation})
grows significantly. \textbf{b)} For $a=20$, the error remains consistently
below $0.3\%$, even as long as $200\; s$ into the simulation. The
tradeoff for this increase in accuracy is a longer computation time.
For both cases the following parameters are used: $\gamma=0.6,\;\alpha=\beta=50\;\frac{units^{2}}{s^{\gamma}},\; dt=0.1s,\;\Delta x=\Delta y=10\; units,\; N_{x}=N_{y}=20.$. We include both maximum error at each timepoint, and average error over the whole solution, at each timepoint.}
\label{fig:a-vs.-accuracy}
\end{center}
\end{figure}

There are various ways to interpret how $a$ is translated into a
numerical algorithm, depending on whether one assumes that it is with
respect to $m$ (mathematical series beginning at 0), time step $n$ (which
may need to be adjusted to begin from index 0 or 1, depending on the
programming language), or time $t$. For example, if one assumes $a$
is a time, we can define $a_{effective}=a/dt$ where $a$ is converted
to a time step. For $dt<1$ this has the same effect as increasing
$a$ and interpreting it with respect to time step $n$. The larger the
value of $a$, the more accurate the scheme presented in Eq. \ref{eq:adaptive-memory},
when compared against the full algorithm in Eq. \ref{eq:2d discretized full equation}.
Fig. \ref{fig:a-vs.-accuracy} shows the error plots where the adaptive
step algorithm in Eq. \ref{eq:adaptive-memory} is compared against
the full implementation (the reference case) Eq. \ref{eq:2d discretized full equation}, for the
adaptive step parameter values $a=4$ and $a=20$. Note that the initial error curve early in the simulation (before the dip and monotonic rise) is a function of $a$.

As with the full implementation, we apply a von Neumann stability
analysis where we assume the solution is separable in the form of Eq. \ref{eq:von-Neumann-F-expansion},
substitute into Eq. \ref{eq:adaptive-memory} and simplify in the same
fashion as described in Section \ref{subsection:vonNeumann} to yield 
\begin{align}
0 = A_{cd}^{n+1}-A_{cd}^{n} +\Delta_{t}^{\gamma}\left( \frac{4\alpha}{\Delta x^{2}}sin^{2}\left(\frac{k_{c}\Delta x}{2}\right)+\frac{4\beta}{\Delta y^{2}}sin^{2}\left(\frac{k_{d}\Delta y}{2}\right)\right)\Bigg\{ \sum_{m=0}^{a}\psi\left(\gamma,m \right) A_{cd}^{n-m} \nonumber \\
 +\sum_{s=2}^{s_{max}(n)}\left( \sum_{\eta=1}^{\eta_{max}(s,n)}\left(2s-1\right)\psi\left(\gamma,M(s,\eta)\right) A_{cd}^{n-M(s,\eta)}+\sum_{p=M(s,\eta_{max})+s}^{min(a^{s},n)}\psi\left(\gamma,p\right)A_{cd}^{n-p}\right) \Bigg\}
\end{align}

As before, we once again assume $A_{cd}^{n+1}=\sigma_{cd}A_{cd}^{n}$
to get the following expression: 
\begin{eqnarray}
\sigma_{cd} & = & 1-\Delta_{t}^{\gamma}\left(\frac{4\alpha}{\Delta x^{2}}sin^{2}\left(\frac{k_{c}\Delta x}{2}\right)+\frac{4\beta}{\Delta y^{2}}sin^{2}\left(\frac{k_{d}\Delta y}{2}\right)\right)\left[\sum_{m=0}^{a}\psi(\gamma,m)\sigma_{uv}^{-m}\right.\nonumber \\
 &  & +\sum_{s=2}^{s_{max}(n)}\left\{ \sum_{\eta=1}^{\eta_{max}\left(s,n\right)}\left(2s-1\right)\psi\left(\gamma,M\left(s,\eta\right)\right)\sigma_{cd}^{-M(s,\eta)}\right.\nonumber \\
 &  & \left.\left.+\sum_{p=M_(s,\eta_{max})+s}^{min(a^{s},n)}\psi(\gamma,p)\sigma_{cd}^{-p}\right\} \right]
\end{eqnarray}

Again, for stability we require that $ $$\mid\sigma_{cd}\mid<1$
for all values of $\sigma$. If we assume the ``worst-case'' scenario
$\left(\sigma_{cd}=-1\right)$, we can infer the following inequality to
provide bounds on the time step and spatial steps 
\begin{eqnarray}
B_{1}(\gamma,n,a) & = & \frac{2}{\Xi}\geq\Delta_{t}^{\gamma}\left(\frac{4\alpha}{\Delta x^{2}}sin^{2}\left(\frac{k_{c}\Delta x}{2}\right)+\frac{4\beta}{\Delta y^{2}}sin^{2}\left(\frac{k_{d}\Delta y}{2}\right)\right)\label{eq:B1-adaptive}\\
\Xi & = & \sum_{m=0}^{a}\psi(\gamma,m)\left(-1\right)^{m}\label{eq:adaptive summation expansion}\\
 &  & +\sum_{s=2}^{s_{max}(n)}\left\{ \sum_{\eta=1}^{\eta_{max}\left(s,n\right)}\left(2s-1\right)\psi\left(\gamma,M(s,\eta)\right)\left(-1\right)^{M(s,\eta)}\right.\nonumber \\
 &  & +\left.\sum_{p=M(s,\eta_{max}+s)}^{min(a^{s},n)}\psi(\gamma,p)\left(-1\right)^{p}\right\} \nonumber 
\end{eqnarray}
If we assume a simple case with a square grid and equal diffusion
coefficients, the most conservative bound occurs when the trigonometric
terms are equal to $1$ and we get 
\begin{equation}
r=\frac{\alpha\Delta t^{\gamma}}{\Delta^{2}}\leq\frac{B_{1}\left(\gamma,n,a\right)}{8}=\frac{1}{4\Xi}=B_{adap}\left(\gamma,n,a\right)\label{eq:Bound-adaptive}
\end{equation}

\begin{figure}
\begin{center}
a)\includegraphics[width=5.8in]{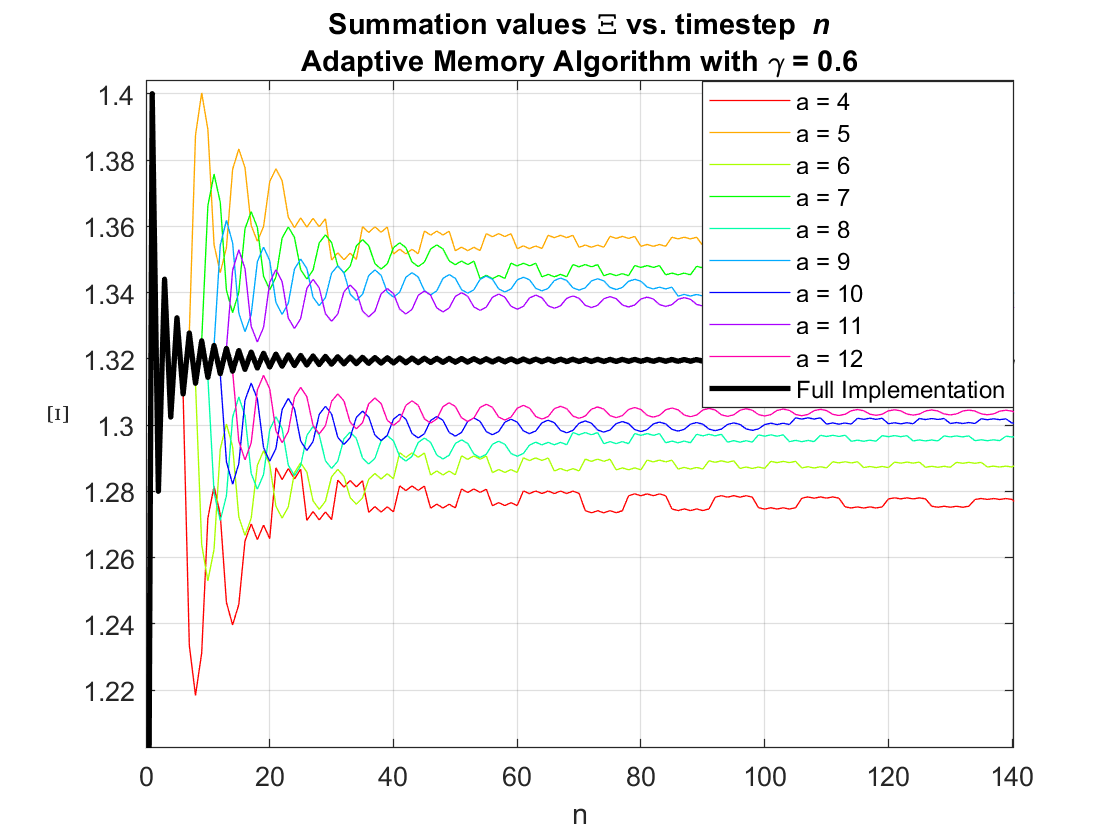} \\ 
\caption[For a range of values for $a$, we plot $\Xi$ as a function of time step $n$. There are
several obvious discontinuities in the pattern of the oscillations,
which occur at the breaks between consecutive $s$ intervals.]

{\footnotesize For a range of values for $a$, we plot $\Xi$ as a function of time step $n$. There are
several obvious discontinuities in the pattern of the oscillations,
which occur at the breaks between consecutive $s$ intervals.}
\label{fig:Xi summation}
\end{center}
\end{figure}

Unlike the full implementation in Section 2.1 where we could use
the convergent nature of the series to approximate the bound $B_{adap}$
with an analytical expression, the adaptive memory algorithm involves
the parameter $a$ which makes it difficult to take a limit approach
to the expression $\Xi$. Fig. \ref{fig:Xi summation}a shows the
value of $\Xi$ for a simulation with parameters $\gamma=0.6,$ and
$a$ varying from $a=4$ to $a=12$. We can see that for each value
of $s$, the interval oscillates around some equilibrium value and
if extended over infinite timesteps, would converge to said value.
As would be expected, we also observe that as $a$ increases, the
value of $\Xi$ approaches the value of the summation in Eq. \ref{eq:B1-full}
in the analogous full 2D case.

From $\Xi$ in Eq. \ref{eq:adaptive summation expansion}, we can take
the sum over $\eta$ for a given value of $s$ 
\begin{equation}
\sum_{\eta=1}^{\eta_{max}\left(s,n\right)}\left(2s-1\right)\psi\left(\gamma,M\left(s,\eta\right)\right)\left(-1\right)^{M\left(s,\eta\right)}\label{eq:alternating-series}
\end{equation}

\noindent and by virtue of the nature of the memory function $\psi$,
\[
\left|\psi\left(\gamma,M_{\eta+1}\right)\right|\leq\left|\psi\left(\gamma,M_{\eta}\right)\right|
\]
which makes Eq. \ref{eq:alternating-series} an alternating series that,
in the limit $\eta_{max}\rightarrow\infty$, converges. However, $\eta$
never approaches infinity in any of these sums, and depending on several of these parameters, the
next $s$ interval can see the sum $\Xi$ oscillating around a different
equilibrium value (see discontinuities in Fig.\ref{fig:Xi summation}).
In many cases it seems that as $n\rightarrow\infty$ and the amplitude
of oscillations becomes smaller and smaller for each subsequent interval,
we do approach a limit or at least a narrow range of values. However,
all of these complications make it difficult to further simplify the
expression for $\Xi$ as an analytical expression.

\begin{figure}
\begin{center}
\includegraphics[width=5.8in]{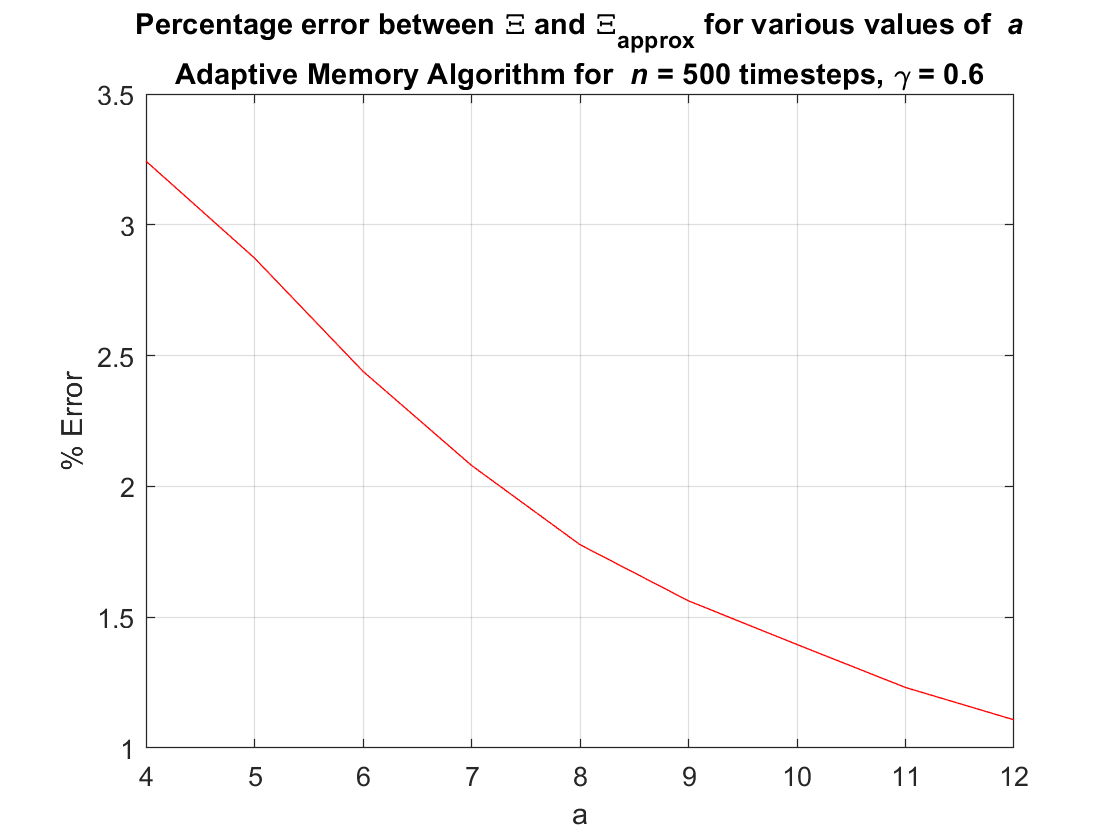} \\ 
\caption[Error values between $\Xi$ and $\Xi_{approx}$ for various values of $a$. These values correspond to $n = 500$ timesteps, and $\gamma =  0.6$.]
{\footnotesize \label{fig:Xi-Xi-approx} Error values between $\Xi$ and $\Xi_{approx}$ for various values of $a$. These values correspond to $n = 500$ timesteps, and $\gamma =  0.6$. }
\end{center}
\end{figure}

On the other hand we do note that because the modulus of the memory
function $\left|\psi\right|$ decreases quickly as $n$ increases,
especially for the first few terms, it is the value of the first summation
term in Eq. \ref{eq:adaptive summation expansion} that largely determines
the final value of $\Xi$ (and in turn $B_{adap}$). Let $\Xi_{approx}=\sum_{m=0}^{a}\psi(\gamma,m)\left(-1\right)^{m}$.
We see from Figure \ref{fig:Xi-Xi-approx} that for various values
of $a$, the error between $\Xi$ (after it approaches a clear equilibrium,
which occurs at least by $n\sim200$ steps in most cases) and $\Xi_{approx}$
is within $2\%$, and as expected, it decreases as $a$ increases. Since
we are unable to extract a simple analytical expression for $\Xi$,
we will instead use $\Xi_{approx}$ which, while numerically based,
is accurate, simple, and quick to calculate for a reasonable range
of $a$ value.

\begin{figure}
\begin{center}
\includegraphics[width=5.8in]{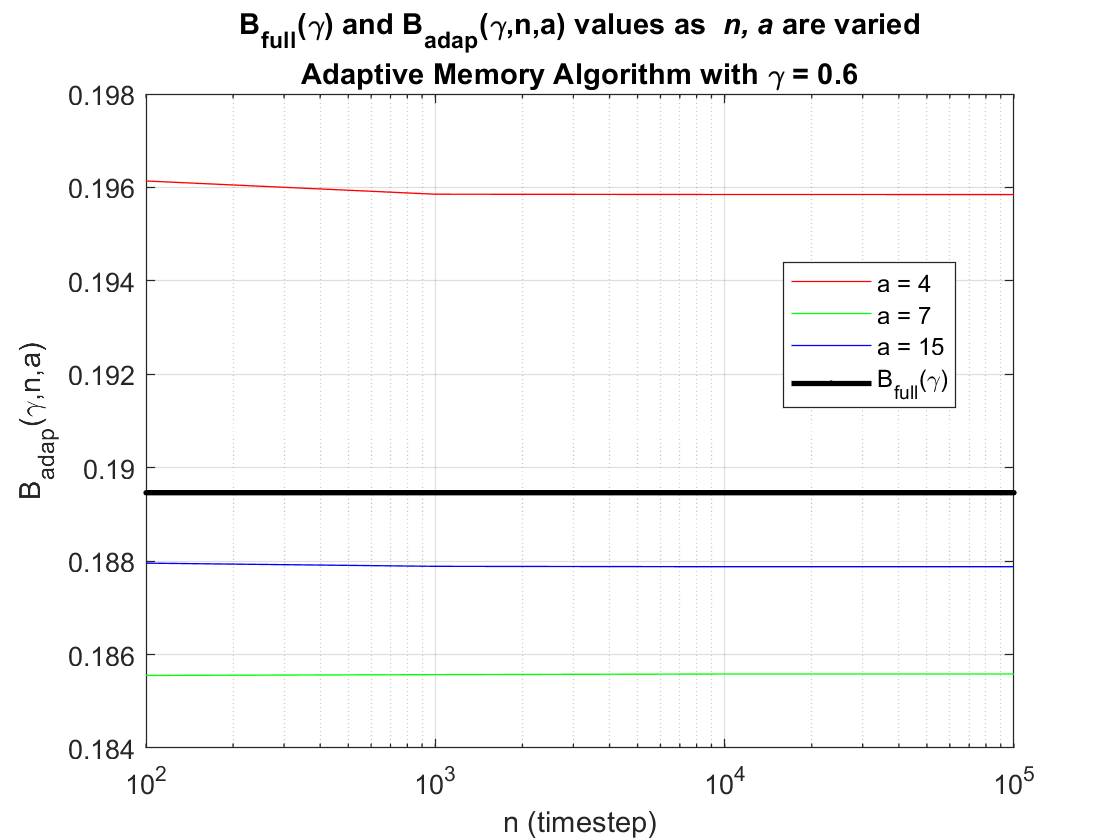} \\
\caption[Adaptive Memory Algorithm with $\gamma = 0.6$: $B_{full}\left(\gamma\right)$ compared to $B_{adap}\left(\gamma,n,a\right)$ values as parameters $a,n$ are varied. The dependence of $B_{adap}$ on $n$ is very weak. As expected, as $a$ increases, $B_{adap}$ approaches $B_{full}$.]
{\footnotesize \label{fig:g=0.6,Bfull,Badap} Adaptive Memory Algorithm with $\gamma = 0.6$: $B_{full}\left(\gamma\right)$ compared to $B_{adap}\left(\gamma,n,a\right)$ values as parameters $a,n$ are varied. The dependence of $B_{adap}$ on $n$ is very weak. As expected, as $a$ increases, $B_{adap}$ approaches $B_{full}$.}
\end{center}
\end{figure}

\begin{figure}
\begin{center}
\includegraphics[width=5.8in]{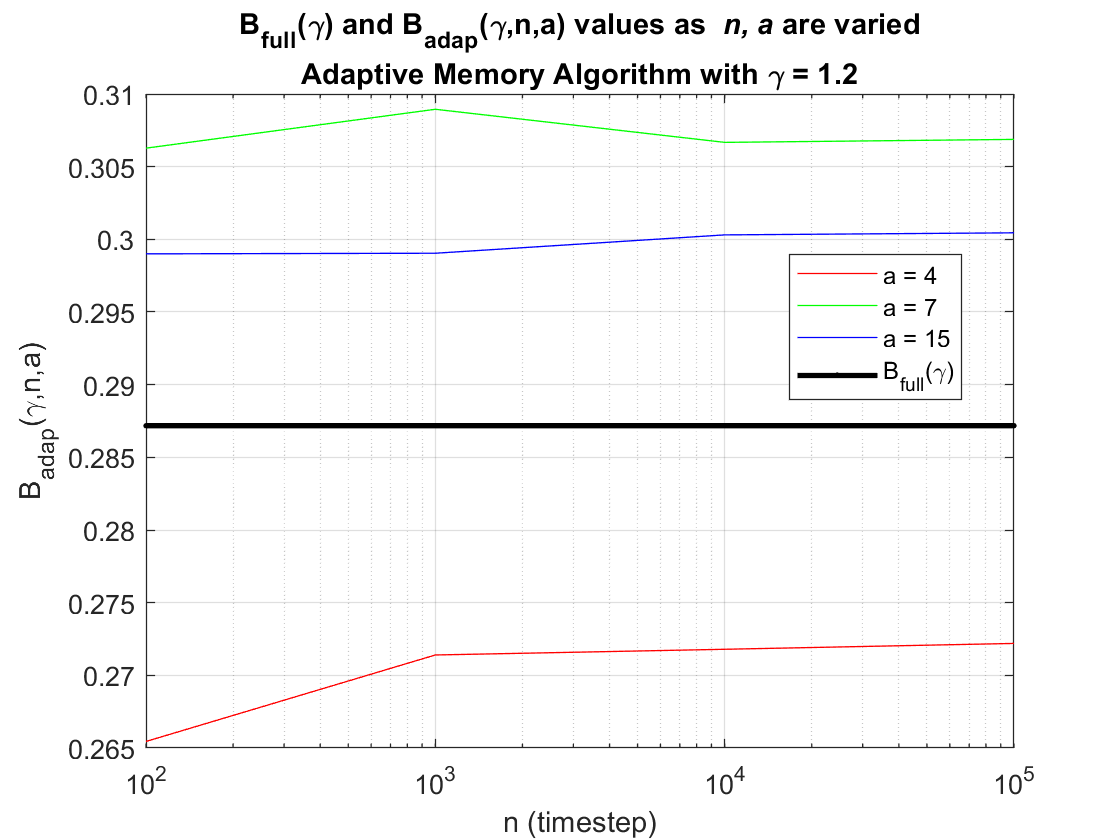} \\
\caption[Adaptive Memory Algorithm with $\gamma = 1.2$: $B_{full}\left(\gamma\right)$ compared to $B_{adap}\left(\gamma,n,a\right)$ values as parameters $a,n$ are varied. $B_{adap}$ a bit more with $n$ compared to the $\gamma = 0.6$ case (Fig. \ref{fig:g=0.6,Bfull,Badap}), but the dependence is still generally weak. As expected, as $a$ increases, $B_{adap}$ approaches $B_{full}$.]
{\footnotesize \label{fig:g=1.2,Bfull,Badap} Adaptive Memory Algorithm with $\gamma = 1.2$: $B_{full}\left(\gamma\right)$ compared to $B_{adap}\left(\gamma,n,a\right)$ values as parameters $a,n$ are varied. $B_{adap}$ a bit more with $n$ compared to the $\gamma = 0.6$ case (Fig. \ref{fig:g=0.6,Bfull,Badap}), but the dependence is still generally weak. As expected, as $a$ increases, $B_{adap}$ approaches $B_{full}$.}
\end{center}
\end{figure}

We will now show with numerical examples how $B_{adap}$ varies with
$\gamma,n,a$, and compare to $B_{full}$, derived in Section \ref{subsection:vonNeumann}.
Figures \ref{fig:g=0.6,Bfull,Badap} and \ref{fig:g=1.2,Bfull,Badap}
show that the values for the $B$ functions change minimally as $n\rightarrow\infty$.
We can also see that $B_{adap}\left(\gamma,n,a\right)$ approaches $B_{full}$ as $a$ increases (which is consistent with Fig. \ref{fig:Xi summation}). As
expected, the actual magnitude of the bounds depend significantly on the order of the fractional operator, $\gamma$. Still, for both the full implementation and the adaptive memory algorithms, it is easy to precompute the $B$ functions based on $\gamma$ (and $a$ in the adaptive case) and then decide how to define
parameters $\Delta t,\Delta x$, and $\Delta y$ as needed to maintain
stability during the simulation. Since the stability dependence on
$n$ is weak in most cases (here we are assuming that time step $n$ is reasonably larger than $a$), the length of simulation has little influence
on the parameters. It is also clear from these results that the adaptive memory bounds $B_{adap}$ are pretty close to $B_{full}$ in many cases, suggesting that the adaptive memory algorithm significantly saves on computational time and power, without adverse effect on the stability regime of the parameters.

\subsection{Simulation Results}

Here we will show that our stability analyses agree well with numerical
simulations. For all plots in the results section, we are using the
adaptive memory algorithm defined in Eq. \ref{eq:adaptive-memory} and
have set $\alpha=50\frac{units^{2}}{s^{\gamma}}$, $\Delta x=\Delta y=10$
$units$, and base interval $a=8$ (which insures a reasonable level
of accuracy with maximum error of the adaptive memory algorithm remaining
under $1\%$ for the duration of the simulation). Our initial condition
is a narrow two-dimensional Gaussian: $u\left(x,y,0\right)=e^{-x^{2}/2\sigma{}_{1}^{2}}e^{-y^{2}/2\sigma_{2}^{2}}$,
$\sigma_{1}=5\; units,\;\sigma_{2}=5\; units$.

We begin by setting $\gamma=0.6$, which puts our simulations in the
subdiffusion regime where the true solution to the fractional diffusion
equation is composed of only decaying modes. Recall that $r=\frac{\alpha\Delta t^{\gamma}}{\Delta^{2}}$, and with the parameters
we have already set, we find that according to Eq. \ref{eq:Bound-adaptive}, for stable conditions, 
$r\leq B_{adap} = 0.1929$ using $\Xi$, and $r\leq B_{adap} = 0.1905$ if we use $\Xi_{approx}$.
Our only degree of freedom is time step $\Delta t$ such that $\Delta t \leq \left( \frac{r\Delta^{2}}{\alpha} \right) ^ {\left( 1/\gamma \right)} = 0.20024 s$ using the more conservative $B_{adap}$ estimate.

\begin{figure}
\begin{center}
a)\includegraphics[height=0.16\paperheight]{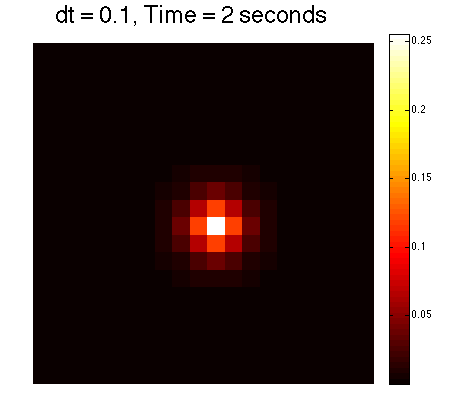}
\includegraphics[height=0.16\paperheight]{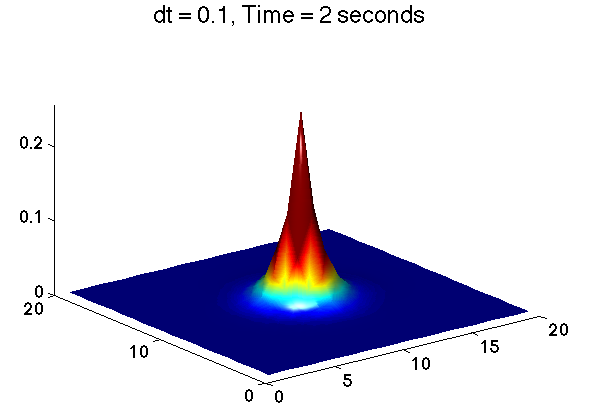}\\

b)\includegraphics[height=0.16\paperheight]{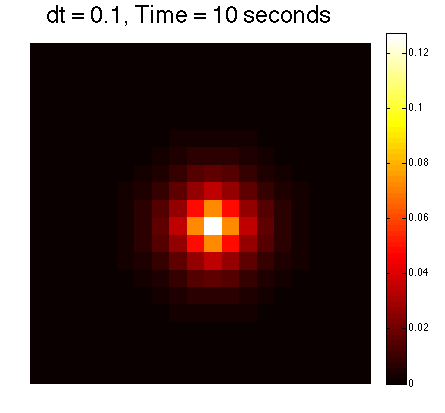}
\includegraphics[height=0.16\paperheight]{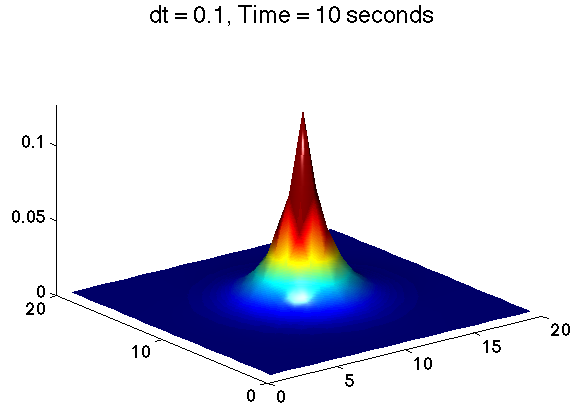}\\

c)\includegraphics[height=0.16\paperheight]{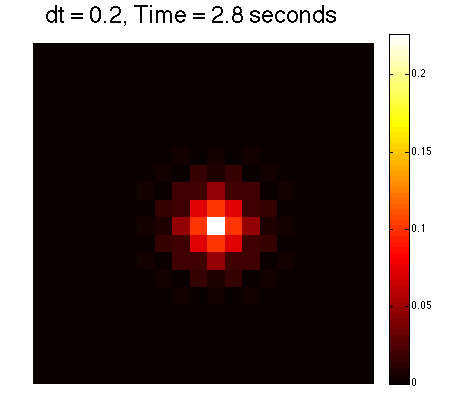}
\includegraphics[height=0.16\paperheight]{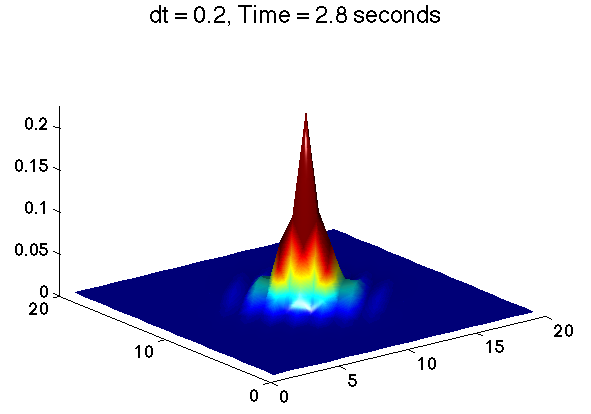}\\

d)\includegraphics[height=0.16\paperheight]{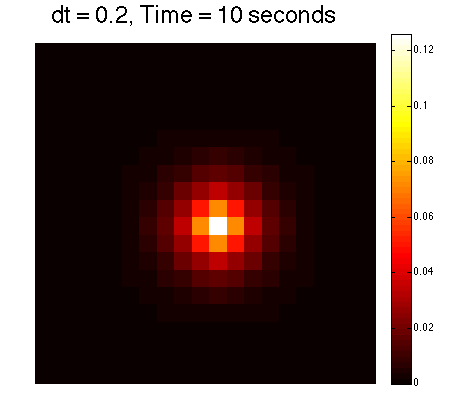}
\includegraphics[height=0.16\paperheight]{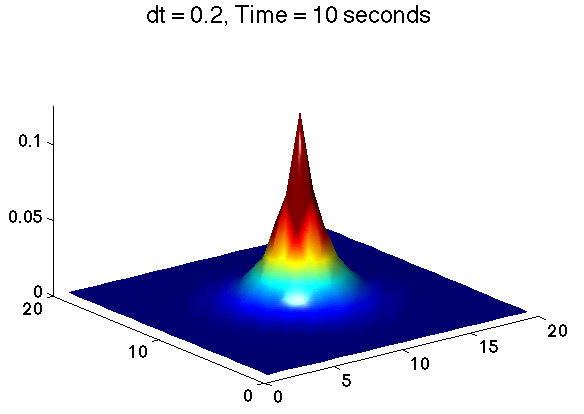} 

\caption[Stable Simulations in the Subdiffusion Region]
{\footnotesize \textbf{Stable Simulations in the Subdiffusion Region}. $\alpha=50\frac{units^{2}}{s^{\gamma}}$,
$\Delta x=\Delta y=10$ $units$, $\gamma=0.6$, left column: 2D
intensity plots, right column: 3D surface plots. \textbf{a,b)} $\Delta t=0.1\; s$,
$r$ is in the stable region and far from the bound, and it is clear
that there is no oscillatory component at any point in the simulation.
\textbf{c,d) $\Delta t=0.2s$}, $r$ in the stable regime, but very
close to the boundary. It is clear in the snapshots, that there
is a tiny oscillatory component that appears early in the simulation
(away from the center of the grid), but decays as the simulation progresses,
and so the solution remains bounded.}
\label{fig:Results-g=0.6-stable}
\end{center}
\end{figure}

In Fig. \ref{fig:Results-g=0.6-stable}$a,b$, we have chosen
$\Delta t=0.1s$ which sets $r = 0.12559$, fulfills the stability inequality criteria and puts our parameters in the stable region of parameter space and relatively
far from the boundary. Both the intensity maps and surface plots confirm
that the numerical solution is clearly stable at early and late time
points in the simulation, as there is no oscillatory behavior.

In Fig. \ref{fig:Results-g=0.6-stable}$c,d$, we have chosen
$\Delta t=0.2s$ ($r = 0.19037$), which puts our parameters in the stable region but close
to the boundary. The plots show that there is an oscillatory component
that is evident early in the simulation. However, at a later time
in the same simulation, the oscillatory component has been suppressed
and the solution has decayed in time and space as expected of subdiffusive
behavior. It is reasonable that the oscillatory component was present
in the beginning of the simulation since $r$ is close to the boundary
between the stable and unstable regimes. However $r$ is still strictly
in the stable regime and this is verified by the fact that the oscillations
do decay with time. Recalling the discussion in Section 2.1.2 relating
to high frequency Fourier modes, we can see clearly from these plots
that the oscillations causing instability are indeed on the order
of the grid discretization.

\begin{figure}
\begin{center}
a)\includegraphics[height=0.19\paperheight]{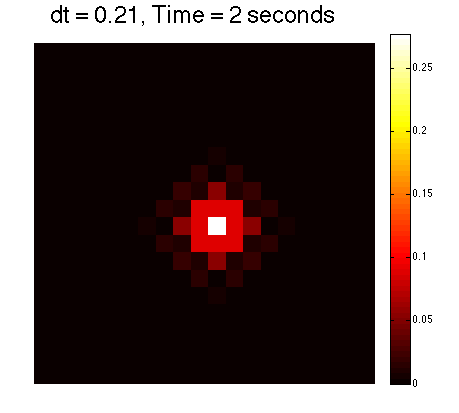}
\includegraphics[height=0.19\paperheight]{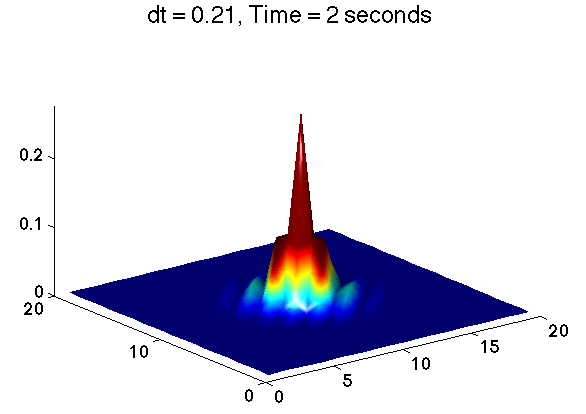}\\

b)\includegraphics[height=0.19\paperheight]{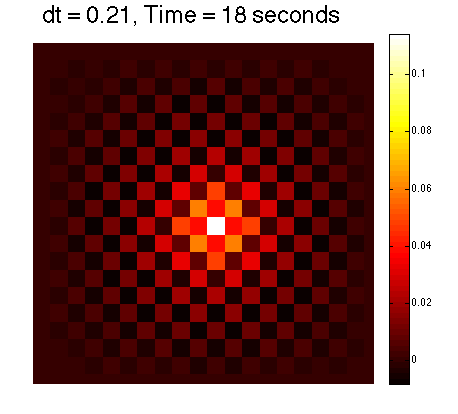}
\includegraphics[height=0.19\paperheight]{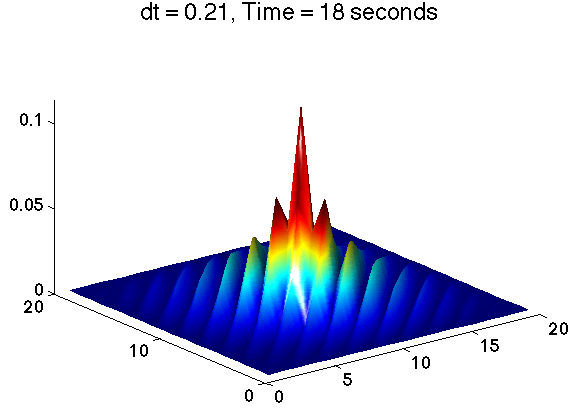}\\

c)\includegraphics[height=0.19\paperheight]{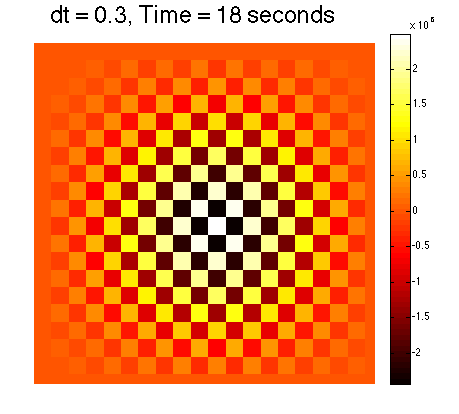}
\includegraphics[height=0.19\paperheight]{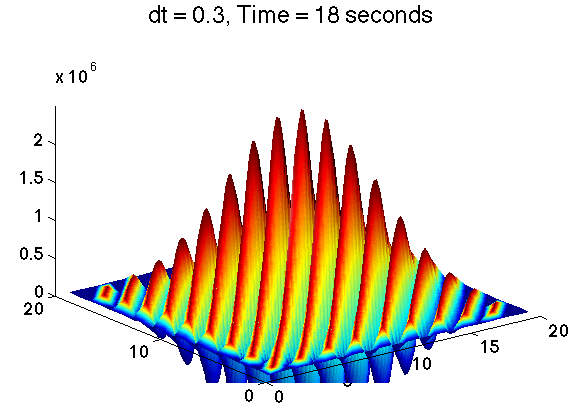}

\caption[Unstable Simulations in the Subdiffusion Region]
{\footnotesize
\textbf{Unstable Simulations in the Subdiffusion Region}.
$\alpha=50\frac{units^{2}}{s^{\gamma}}$,
$\Delta x=\Delta y=10$ $units$, $\gamma=0.6$, left column: 2D intensity
plots, right column: 3D surface plots. \textbf{a,b)} $\Delta t=0.21s$, $r$
is in the unstable region and close to the boundary. The plots reflect
an oscillatory component that appears early in the simulation and
is sustained through the duration of the simulation. Even though the
oscillations do not increase, and do not overwhelm the decaying modes
of the true solution, they do not disappear and would remain even
as $n\rightarrow\infty$, which is an indication of instability. \textbf{c)
$\Delta t=0.3s$}, $r$ is in the unstable region, and far from the
bounds between the unstable and stable regions. The plots reflects
oscillatory components that remain at a later time in the simulation and have increased in an unbounded manner. If we note the axis scale of both plots, we also observe that those components
have completely overwhelmed the decaying modes of the true solution. These are a hallmarks
of strong instability in numerical simulations.}
\label{fig:Results-g=0.6-unstable}
\end{center}
\end{figure}

In Fig. \ref{fig:Results-g=0.6-unstable}$a,b$, we have chosen $\Delta t=0.21s$, which sets $r = 0.19602$ and puts our parameters in the unstable region,
but close to the boundary between stable and unstable regions. As expected of  unstable simulations, the
numerical solution includes oscillatory components that appear early in the simulation and persist as the simulation
continues (they would remain even as $n\rightarrow\infty$). But because $r$ is close to the boundary with the stable
regime, the oscillations do not drastically overwhelm
the decaying modes, suggesting that the solution is only mildly unstable.

In Fig. \ref{fig:Results-g=0.6-unstable}$c$, we have
chosen $\Delta t=0.3s$, which sets $r = 0.2428$, putting our parameters clearly in the unstable
region and further from the boundary. The plots depicting a later time in the simulation confirm this by showing
that the oscillatory behavior of the solution grows unbounded with
time and completely overwhelms the decaying modes of the true solution,
as observed by scale of the $z$ axis, which is orders of magnitude larger than
that of the true solution.

Next, we will consider the superdiffusion cases with $\gamma=1.2$.
According to Eq. \ref{eq:Bound-adaptive}, for stable numerical solutions,
$r\leq B_{adap} = 0.272$ using $\Xi$, and $r\leq B_{adap} = 0.2809$ if we use $\Xi_{approx}$.
In addition, it is less straightforward to characterize the superdiffusion regime because in addition to decaying modes characteristic of diffusion,
the true solution also has oscillatory components (which is logical because
as $\gamma$ increases towards 2, we approach the classical wave equation). However,
we can still observe the difference in behavior between unstable and
stable numerical solutions.

\begin{figure}
\begin{center}
a)\includegraphics[height=0.16\paperheight]{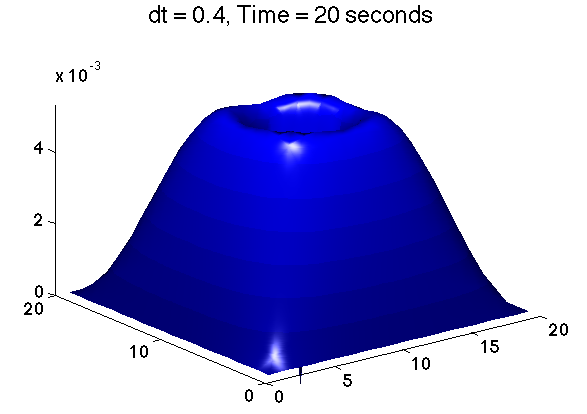}
b)\includegraphics[height=0.16\paperheight]{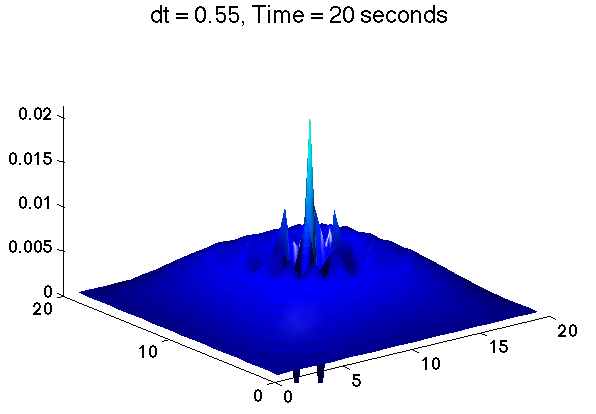}\\
c)\includegraphics[height=0.16\paperheight]{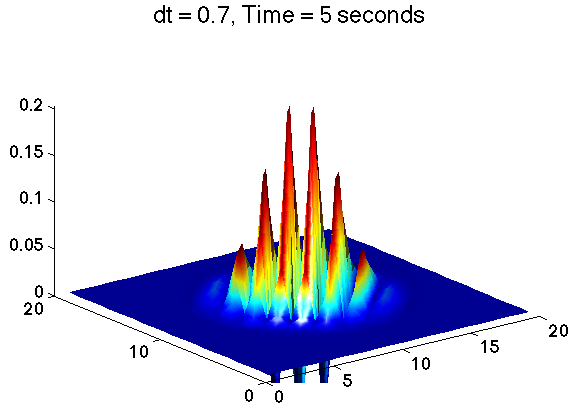}
d)\includegraphics[height=0.16\paperheight]{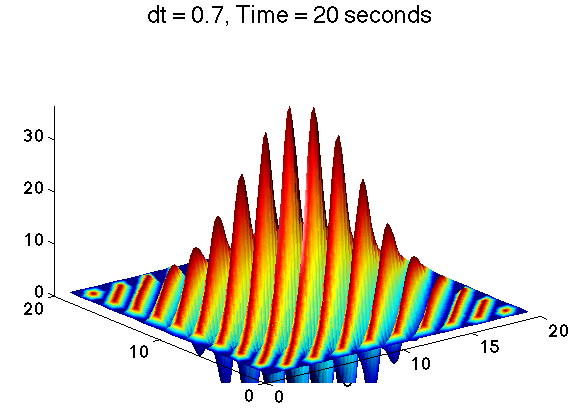}
\caption[Stable and Unstable Superdiffusion Simulations]{\footnotesize

\textbf{Stable and Unstable Superdiffusion Simulations}. $\alpha=50\frac{units^{2}}{s^{\gamma}}$,
$\Delta x=\Delta y=10$ $units$, $\gamma=1.2$. \textbf{a)} $\Delta t=0.4s$, $r$
is in the stable regime. Because we are in the superdiffusion regime,
the true solution does have oscillatory components alongside decaying
modes. However, the numerical solution is obviously bounded, even
later in the simulation, and there is no sign of high frequency oscillations
(visible on the plots \textbf{c} and \textbf{d}), so the plots seem to verify that the
simulation is indeed stable. \textbf{b)} $\Delta t=0.55s$, $r$ is
in the stable regime. We can clearly see high frequency components,
especially contrasted the same plot at the same time in the simulation,
in part \textbf{a}. However, as observed by the scale bar, the high frequency
components don't overwhelm the decaying modes of the solution and
cause the solution to grow unbounded, so this verifies stability. \textbf{c,d)
$\Delta t=0.7s$},\textbf{ $r$ }is in the unstable regime. It is
clear from the plots that the high frequency oscillatory components
is present early in the simulation, but the decaying modes still dominate
the overall solution. However, it's clear that later in the simulation
the high frequency oscillations continue to grow unbounded in time,
which verifies the unstable nature of the numerical solution.} 
\label{fig:Results-g=1.2}
\end{center}
\end{figure}

In Fig. \ref{fig:Results-g=1.2}$a$, we have chosen $\Delta t=0.4s$,
which sets $r = 0.16651$, putting our parameters in the stable region. The plots show a numerical solution
with low spatial frequency oscillatory components (oscillating on a much larger scale than the scale
of the grid elements), and the presence of high frequency components characteristic
of noise, is absent. The numerical solution remains bounded. In Fig. \ref{fig:Results-g=1.2}$b$, $\Delta t=.55s$,
which sets $r = 0.24401$, putting our parameters still in the stable region but close to the boundary between stable and unstable regions.
The presence of the high frequency noise components is apparent, especially compared to Fig. \ref{fig:Results-g=1.2}$a$ at the same simulation time,  but
even so late into the simulation, they do not overwhelm the overall solution
and are clearly bounded, which verifies stability. Finally, in Fig. \ref{fig:Results-g=1.2}$c,d$,
we have chosen $\Delta t=0.7s$, setting $r = 0.3259$ and putting our parameters clearly in the unstable
region. The plots show high frequency oscillatory behavior that appears early in the simulation when decaying modes still dominate the solution. At a later timestamp the oscillatory components have grown drastically and in an unbounded manner, which,
as also noted in Fig. \ref{fig:Results-g=0.6-unstable}, is
a good indicator of numerical instability.

All our results verify the conclusions made during our stability analysis
in Sections 2.1.2 and 2.2. It's also clear from Figs. \ref{fig:Results-g=0.6-stable}
and \ref{fig:Results-g=0.6-unstable} that the boundary
between stable and unstable regimes is quite a sharp one, as changing
our time step from $\Delta t=0.2s$ to $\Delta t=0.21s$ was enough
to change our simulations from being stable to unstable.

%

%

\section{Complexity Analysis}

In the following sections we analyze the complexity of the full two-dimensional
implementation of the fractional diffusion equation, the adaptive
memory algorithm, and the linked list alternative version defined
in \cite{Paper0}.

\subsection{Full Implementation}

We begin with the pseudocode for the full implementation defined in Eq. \ref{eq:2d discretized full equation}. In our pseudocode we generally omit instructions not essential to showing the logical structure of the algorithm and note that
the most important aspect of time complexity analysis are data structures
involving loops. Lines of code that are run a constant number of times and are independent of variables like timestep, can be ignored in the overall Big-$O$
functional form. 

From the pseudocode in Algorithm \ref{full-pseudocode} we can see that
for each $n$, the total approximate number of instructions (ignoring
things like constant number of instructions) is $N_{x}N_{y}\left(n+1\right)$ where
$n$ varies from $0:N$ . Therefore:
\begin{align*}
total\; instructions & =N_{x}N_{y}\left(1+...+N+1\right)\\
 & =N_{x}N_{y}\left(\frac{\left(N+1\right)\left(N+2\right)}{2}\right)\\
 & \Rightarrow O\left(N_{x}N_{y}N^{2}\right)
\end{align*}

This gives us the growth in execution time as a function of the grid
dimensions $N_{x}$ and $N_{y}$, as well as the number of timesteps $N$. While
execution time grows linearly with the size of each spatial dimension of the
grid we are using, it grows as the square of the number of timesteps
in the simulation, which is not the most efficient numerical implementation.

\begin{algorithm}
	\caption{Full 2D Implementation}\label{full-pseudocode}
	\begin{algorithmic}[1]
		\For {$n = 0:N$} \Comment{$N$ is the total \# of timesteps the algorithm will run}
			\For{$j = 0:N_x-1$} \Comment{iterate over all grid points in x direction}
				\For{$l = 0:N_y-1$} \Comment{iterate over all grid points in y direction}
					\For{$m = 0:1:n$} 
						\State Calculate $\psi(\gamma,m)\left(\frac{\alpha}{\Delta x^{2}}\delta x_{j,l}^{n-m}+\frac{\beta}{\Delta y^{2}}\delta y_{j,l}^{n-m}\right)$ 
					\EndFor 
					\State Calculate $u_{j,l}^{n+1}$ 
				\EndFor
			\EndFor
		\EndFor 
	\end{algorithmic}
\end{algorithm}

Here we note that in the one-dimensional case, we can easily vectorize the loop iterating over grid points, increasing efficiency with regards to grid size. However, with respect to the time dimension, complexity still grows as $O\left(N^{2}\right)$.
\subsection{Adaptive Memory Algorithm}

In Algorithm 2 we provide the pseudocode for the implementation of the adaptive
memory algorithm described in Eq. \ref{eq:adaptive-memory}. As before,
we omit the instructions detailing the calculations themselves because
counting these instructions involve multiplying constants against
the number of times the loops are run, and do not affect the form
of the Big-$O$ expressions in terms of the input variable
of interest ($N_x, N_y, N$). As with the full implementation, we briefly note that in the one-dimensional analog, we can vectorize the loop iterating over spatial gridpoints.

Analysis of the total number of instructions in Part 1 of Algorithm \ref{am-pseudocode}
proceeds in exactly the same way as in Section 3.1. The total instruction
count is approximately
\begin{align*}
N_{x}N_{y}\left(1+2+...+a+1\right) & =N_{x}N_{y}\frac{\left(a+1\right)\left(a+2\right)}{2}\\
 & \Rightarrow O\left(N_{x}N_{y}a^{2}\right)
\end{align*}

\begin{algorithm}
	\caption{Adaptive Memory Algorithm}\label{am-pseudocode}
	\begin{algorithmic}[1]
		\State \Comment{\textbf{Part 1}}
		\For{$n = 0:a$} 
			\For{$j = 0:N_x-1$} \Comment{iterate over all grid points in x direction}
				\For{$l = 0:N_y-1$} \Comment{iterate over all grid points in y direction}
					\For{$m = 0:1:n$} 
						\State Calculate $\psi(\gamma,m)\left(\frac{\alpha}{\Delta x^{2}}\delta x_{j,l}^{n-m}+\frac{\beta}{\Delta y^{2}}\delta y_{j,l}^{n-m}\right)$
					\EndFor 
					\State At this location, calculate $u_{j,l}^{n+1}$ 
				\EndFor
			\EndFor
		\EndFor
		
		\State \Comment{\textbf{Part 2}}
		\For {$n = a+1:N$} \Comment{$N$ is the total \# of timesteps the algorithm will run} \tikzmark{top3} \tikzmark{right3}
			\State Determine $s_{max}$ based on $n$
			\For{$j = 0:N_x-1$} \Comment{iterate over all grid points in x direction}
				\For{$l = 0:N_y-1$} \Comment{iterate over all grid points in y direction}
				
					\For{$m = 0:1:a$} \Comment{Calculate Sum 1 (base interval)}
						\State Calculate R1: $\psi(\gamma,m)\left(\frac{\alpha}{\Delta x^{2}}\delta x_{j,l}^{n-m}+\frac{\beta}{\Delta y^{2}}\delta y_{j,l}^{n-m}\right)$ \Comment{R1 from Eq. \ref{eq:adaptive-memory}}
					\EndFor 
					
					\For{$s = 2:s_{max}$} 
						\For{$\eta = 1:\eta_{max}$} \tikzmark{top1}
							\State Calculate R2 in Eq. \ref{eq:adaptive-memory}\tikzmark{right1} 
						\EndFor \tikzmark{bottom1}
						
						\For{$p = M_{max}+s:min\left(a^s,n\right)$} \tikzmark{top2}\tikzmark{right2}
							\State Calculate R3 in Eq. \ref{eq:adaptive-memory}
						\EndFor \tikzmark{bottom2}
					\EndFor
					\State Use R1-R3 to calculate $u_{j,l}^{n+1}$ 
				\EndFor
			\EndFor
		\EndFor \tikzmark{bottom3}
	\end{algorithmic}
	\AddNote{top1}{bottom1}{right1}{ $\{L1\}$ }
	\AddNote{top2}{bottom2}{right2}{ $\{L2\}$ }
	\AddNoteLeft{top3}{bottom3}{right3}{ $\{L3\}$ }
\end{algorithm}

Analysis of Part 2 of Algorithm \ref{am-pseudocode} is more complicated
because of the additional parameters $s_{max},\;\eta_{max}, $ and $ M_{max}$.
For loop \{L1\} (lines 21-23), $\eta_{max}$ depends on the current $s$ interval, but the maximum value is always $\frac{a^{s}-a^{s-1}}{2s-1}$.
For loop \{L2\} (lines 24-26), the loop serves to sample timesteps when the current
interval $[a^{s-1}:a^{s}]$ (or $[a^{s-1}:n]$ for the latest interval corresponding to $s_{max}$), is not evenly divisible by the weight for
that interval, $(2s-1)$, and therefore there are a few timesteps that aren't
taken into account with loop \{L1\}. The maximum number of times
loop \{L2\} is run, is always $<2s-1$. Therefore in the large \textit{for} loop \{L3\}
in Algorithm \ref{am-pseudocode} encompassing lines 13-31, for a given value of $n$, the number
of instructions run is, at most: 
\begin{equation}
total\;instructions \leq N_{x}N_{y}\left\{ a+1 +\sum_{s=2}^{s_{max}}\left(\underbrace{\frac{a^{s}-a^{s-1}}{2s-1}}_{\{1a\}}+\underbrace{2s-1}_{\{1b\}}\right)\right\} \label{eq:total-instructions-am-1}
\end{equation}
To simplify this expression further we need to make some assumptions about worst
case scenarios. We can simplify the first term \{1a\} of the summation in Eq. \ref{eq:total-instructions-am-1} using telescoping series:
\[
\]
\begin{align*}
\sum_{s=2}^{s_{max}}\frac{a^{s}-a^{s-1}}{2s-1} & =\frac{a^{2}-a}{3}+\frac{a^{3}-a^{2}}{5}+...+\frac{a^{s_{max}}-a^{s_{max}-1}}{2s_{max}-1}\\
 & <\frac{a^{2}-a}{3}+\frac{a^{3}-a^{2}}{3}+...+\frac{a^{s_{max}}-a^{s_{max}-1}}{3}\\
 & =\frac{1}{3}\left\{ \left(a^{2}-a\right)+\left(a^{3}-a^{2}\right)+...+\left(a^{s_{max}-1}-a^{s_{max}-2}\right)+\left(a^{s_{max}}-a^{s_{max}-1}\right)\right\} \\
 & =\frac{1}{3}\left(a^{s_{max}}-a\right)
\end{align*}
The second term \{1b\} of the summation can be reduced to
\begin{align*}
\sum_{s=2}^{s_{max}}2s-1 & =2\sum_{s=2}^{s_{max}}s-\sum_{s=2}^{s_{max}}1\\
 & =2\left(\frac{s_{max}\left(s_{max}+1\right)}{2}-1\right)-\left(s_{max}-1\right)\\
 & =s_{max}^{2}-1
\end{align*}
The expression in Eq. \ref{eq:total-instructions-am-1} can now be reduced
to
\begin{equation}
total\;instructions \leq N_{x}N_{y}\left( a+\frac{1}{3}a^{s_{max}}-\frac{1}{3}a+s_{max}^{2}-1\right) \label{eq:total-instructions-am-2}
\end{equation}
$s_{max}$ for a given timestep $n$ is determined by $a^{s_{max}-1}+1\leq n\leq a^{s_{max}}$, or $\frac{ln\left(n\right)}{ln\left(a\right)}\leq s_{max}\leq\frac{ln(n-1)}{ln(a)}+1$.
Since we are concerned with worst case scenarios, we use the maximum
value of this interval, $s_{max}=\frac{ln(n-1)}{ln(a)}+1$ and substitute
into Eq. \ref{eq:total-instructions-am-2} . Now we must consider the
\textit{for} loop \{L3\} in Algorithm 2 and how the total number of instructions relates to
the total number of timesteps $N$: 
\begin{equation}
total\;instructions \leq N_{x}N_{y}\sum_{n=a+1}^{N}\left( \underbrace{\frac{2}{3}a}_{\{2a\}}+\underbrace{\frac{1}{3}a^{ln(n-1)/ln(a)}}_{\{2b\}}a+\underbrace{\left(\frac{ln(n-1)}{ln(a)}\right)^{2}}_{\{2c\}}+\underbrace{2\frac{ln(n-1)}{ln(a)}}_{\{2d\}}\right) \label{eq:total-instructions-am-3}
\end{equation}
We now simplify the four components of the summation in 
Eq. \ref{eq:total-instructions-am-3}. 

For \{2a\}, $\sum_{n=a+1}^{N}\frac{2}{3}a=\left(N-a\right)\frac{2}{3}a$. Since $a$ is a constant parameter chosen by the user prior to the
simulation, the dominant term in this expression is clearly $N$.

For \{2b\}:
\begin{align*}
\frac{a}{3}\sum_{n=a+1}^{N}a^{\frac{ln(n-1)}{ln(a)}} & =\frac{a}{3}\sum_{n=a+1}^{N}a^{\frac{log_{a}\left(n-1\right)}{log_{a}e}\frac{1}{lna}}\\
 & =\frac{a}{3}\sum_{n=a+1}^{N}\left[a^{log_{a}\left(n-1\right)}\right] ^ \frac{1}{\left(log_{a}e\right)\left( log_{e}a\right)}\\
 & =\frac{a}{3}\sum_{n=a+1}^{N}n-1\\
 & =\frac{a}{3}\left(\frac{N^{2}+N}{2}-\left(1+2+...+a\right)-\left(N-a\right)\right)
\end{align*}
The dominant term here is clearly $N^{2}$.

For \{2c\} in equation Eq. \ref{eq:total-instructions-am-3}:
\begin{align*}
\frac{1}{\left(ln(a)\right)^{2}}\sum_{n=a+1}^{N}\left(ln(n-1)\right)^{2} & =\frac{1}{\left(ln(a)\right)^{2}}\left[\left(ln(a)\right)^{2}+\left(ln(a+1\right)^{2}+...+\left(ln(N-1)\right)^{2}\right]\\
 & \leq\frac{1}{\left(ln(a)\right)^{2}}\left[\left(ln(N-1)\right)^{2}+\left(ln(N-1)\right)^{2}+...+\left(ln(N-1)\right)^{2}\right]\\
 & =\frac{1}{\left(ln(a)\right)^{2}}\left[\left(N-a\right)\left(ln(N-1)\right)^{2}\right]
\end{align*}
This expression in Big-$O$ notation is $O\left(N\left(ln\left(N\right))^{2}\right)\right)$.
Alternatively, we can simplify \{2c\} using integration by parts, to
get the same result. 

For \{2d\}:
\begin{align*}
\frac{2}{ln(a)}\sum_{n=a+1}^{N}ln(n-1) & =\frac{2}{ln(a)}ln\left[\left(a\right)\left(a+1\right)...\left(N-1\right)\right]\\
 & =\frac{2}{ln(a)}ln\left(\frac{\left(N-1\right)!}{\left(a-1\right)!}\right)\\
 & =\frac{2}{ln(a)}\left(ln\left(\left(N-1\right)!\right)-ln\left(\left(a-1\right)!\right)\right)
\end{align*}
Using Stirling's approximation, our final expression is 
\begin{align}
\frac{2}{ln(a)}\left(Nln(N-1)-N+1+O\left(ln\left(N-1\right)\right)-ln(\left(a-1\right)!)\right)
\end{align}
which is dominated by the $Nln\left(N-1\right)$term. 

Combining the simplified expressions for \{2a\}-\{2d\} in Eq. \ref{eq:total-instructions-am-3},
it is clear that the $N^{2}$ term dominates the behavior of the
overall expression as $N$ becomes large, and so we conclude that
the worst case behavior of the adaptive step algorithm, is still bounded
by $O\left(N_{x}N_{y}N^{2}\right)$. While the adaptive
time step algorithm improves raw execution time in relation to the
full memory implementation, the complexity analysis shows that when
it comes to scaling with large $N$, the algorithm is as inefficient
as the full implementation.

\subsection{Linked List Adaptive Timestep Algorithm}

In \cite{Paper0} we describe an alternative version of the adaptive memory
algorithm which is based on a power law that enables us to eliminate
past history timepoints that will never again need to be referenced,
thus saving us execution time and memory. Here we summarize the algorithm and
refer to \cite{Paper0, Thesis} for additional background technical details
and derivation. 
\begin{align}
\frac{u_{j,l}^{n+1}-u_{j,l}^{n}}{\Delta_{t}^{\gamma}} & =\sum_{\left\{ u_{j,l}^{i}\in U^{n}\right\} }\psi(\gamma,n-i)w^{i}\left(\frac{\alpha}{\Delta x^{2}}\delta x_{j,l}^{i}+\frac{\beta}{\Delta y^{2}}\delta y_{j,l}^{i}\right)\label{eq:linked-list-def}\\
w^{n+1} & =1\nonumber \\
W^{n+1} & :=\left\{ w^{i}\in W^{n}\right\} +\left\{ w^{n+1}\right\} \nonumber \\
U^{n+1} & :=\left\{ u^{i}\in U^{n}\right\} +\left\{ u^{n+1}\right\} \nonumber 
\label{ll-def}
\end{align}
The data associated with each timestep $n$ is stored in the
elements constituting a doubly linked list, and $i$ represents the timestep in each node/element in the linked list. Additionally,
when there are more than $\eta$ points in the set $W^{i}$ of any
given weight, the elements of the set are condensed according to Algorithm
4.1 in \cite{Paper0}. We present the pseudocode for the linked
list implementation in Algorithm \ref{ll-pseudocode}.

\begin{algorithm}
	\caption{Linked List Implementation}\label{ll-pseudocode}
	\begin{algorithmic}[1]
		\For {$n = 0:N$} \Comment{$N$ is the total \# of timesteps the algorithm will run}
			\For{$j = 0:N_x-1$} \Comment{iterate over all grid points in x direction}
				\For{$l = 0:N_y-1$} \Comment{iterate over all grid points in y direction}
					\For{$node = 1$:length(linked list)}  \tikzmark{top1}
						\State From each node calculate $\psi(\gamma,n-i)w^i\left(\frac{\alpha}{\Delta x^{2}}\delta x_{j,l}^{i}+\frac{\beta}{\Delta y^{2}}\delta y_{j,l}^{i}\right)$ \tikzmark{right1}
					\EndFor \tikzmark{bottom1}
					\State Calculate $u_{j,l}^{n+1}$ 
				\EndFor
			\EndFor

			\State Create a new node in the list to store $u_{j,l}^{n+1}$
			\State If-else statements to determine whether we need to condense the linked list according to Algorithm 4.1 in \cite{Paper0}. If so:
			\While{it continues to be necessary to condense the linked list} \tikzmark{top2} 
				\State Constant number of instructions to delete nodes and update weights \tikzmark{right2}
			\EndWhile \tikzmark{bottom2}
		\EndFor
	\end{algorithmic}
	\AddNote{top1}{bottom1}{right1}{ $\{L1\}$ }
	\AddNote{top2}{bottom2}{right2}{ $\{L2\}$ }
\end{algorithm}

To determine how many times loop \{L1\} runs in relation to timestep
$n$, we can consider the following. For a given timestep $n$, we
are interested in a summation of terms involving weights (powers of
two), each set of which never exceeds $\eta$ nodes. For example,
for timestep $n=25$ and $\eta=5$, Table \ref{tab:Linked-list-example} shows the timestep and weight for all nodes currently in the linked list at this point in the simulation. 
\begin{table}
\noindent \centering{}\caption[Linked list data for timestep $n=25$]{\footnotesize Linked list data for timestep $n=25$\label{tab:Linked-list-example}}
\begin{tabular}{|c|c|c|c|c|c|c|c|c|c|c|c|c|}
\hline 
Timestep & $0$ & $4$ & $8$ & $12$ & $14$ & $16$ & $18$ & $20$ & $22$ & $23$ & $24$ & $25$\tabularnewline
\hline 
\hline 
Weight & $2^{2}=4$  & $4$ & $4$ & $2^{1}=2$ & $2$ & $2$ & $2$ & $2$ & $2^{0}=1$ & $1$ & $1$ & $1$\tabularnewline
\hline 
\end{tabular}
\end{table}

Since the dynamically changing weight is directly related to how often
the history of the simulation is sampled, we can construct the following
relationship between time step and weight based on the data in Table \ref{tab:Linked-list-example}:
\[
2^{0}\cdot{5}+2^{1}\cdot{5}\leq n+1=26\leq2^{0}\cdot{5}+2^{1}\cdot{5}+2^{2}\cdot{5}
\]
The inequality arrives from the fact that a particular weight category
has $\eta$ or fewer nodes. 

For a general timestep $n$ we can write the inequalities in terms
of geometric series:
\begin{equation}
\sum_{l=0}^{x-1}2^{l}\leq\frac{n+1}{\eta}\leq\sum_{l=0}^{x}2^{l}\label{eq:ll-inequality}
\end{equation}
where $x+1$ is the number of weight categories represented in the
linked list at timestep $n$. In the example in Table \ref{tab:Linked-list-example},
the number of weight categories is 3, with weights $2^{0},\;2^{1},\;\mbox{and}\;2^{2}$.
Substituting the geometric series with their total sums and taking $log_{2}$ of both sides, we can rewrite
equation Eq. \ref{eq:ll-inequality} as 
\begin{align}
x & \leq log_{2}\left(\frac{n+1}{\eta}+1\right)\leq x+1\Rightarrow\nonumber \\
log_{2}\left(\frac{n+1}{\eta}+1\right) & \leq x+1\leq log_{2}\left(\frac{n+1}{\eta}+1\right)+1\label{eq:max-weight-categories}
\end{align}
where $x+1$ is the number of weight categories. 

We now return to Algorithm \ref{ll-pseudocode} and refer to loop
\{L1\}. The maximum possible number of times this loop runs for any
given timestep $n$, is the length of the linked list, which is $\eta$ times the maximum number of weight
categories at that time. Using Eq. \ref{eq:max-weight-categories}, we get
\[
\eta\left(log_{2}\left(\frac{n+1}{\eta}+1\right)+1\right)
\]
Therefore the total number of instructions in lines 2-9 of Algorithm \ref{ll-pseudocode} is 
$N_{x}N_{y}\eta\left(log_{2}\left(\frac{n+1}{\eta}+1\right)+1\right)$.
We now look at loop \{L2\} (lines 11-13). This \textit{while} loop runs whenever the number of nodes belonging to
a particular weight category, exceeds the user-set parameter $\eta$, and we need to condense the linked list.
The maximum number of times this may run within a single timestep,
is the number of weight categories in the linked list. For example,
if the number of nodes with weight $2^{0}$ exceeds $\eta$, we delete
the second ``least'' node in the set of nodes with this weight,
and double the weight of the ``least'' element (node) in this same
set, to $2^{1}$. (See the derivation of Algorithm 4.1 in \cite{Paper0} for details about
the ordering of these sets). For some particular timestep this may
mean that the number of nodes with weight $2^{1}$ now exceeds $\eta$
and we continue to condense in the same fashion until we have addressed
all weight categories currently represented in the linked list. The
maximum number of weight categories is given by Eq. \ref{eq:max-weight-categories},
and this is the maximum number of times loop \{L2\} is run for a given time step
$n$. Now we can combine all loops in Algorithm \ref{ll-pseudocode} and
write that the maximum number of instructions (ignoring constants)
is given by
\begin{align}
total\; instructions & \leq \sum_{n=0}^{N}\left[N_{x}N_{y}\eta\left(log_{2}\left(\frac{n+1}{\eta}+1\right)+1\right)+log_{2}\left(\frac{n+1}{\eta}+1\right)+1\right]\nonumber \\
& =\left(N_{x}N_{y}\eta+1\right)\left[\underbrace{\sum_{n=0}^{N}log_{2}\left(\frac{n+1}{\eta}+1\right)}_{\{1\}}+N+1\right]\label{eq:ll-total-instructions}
\end{align}
We can simplify the summation term \{1\}:
\[
\sum_{n=0}^{N}log_{2}\left(\frac{n+1+\eta}{\eta}\right)=\sum_{n=0}^{N}\left(log_{2}\left(n+1+\eta\right)-log_{2}\left(\eta\right)\right)
\]
With a change of variable $d=n+1+\eta$ we can rewrite the summation and once again utilize Sterling's approximation to arrive at:
\begin{align*}
\sum_{d=1+\eta}^{N+1+\eta=N'}log_{2}\left(d\right)-Nlog_{2}\left(\eta\right) & \leq log_{2}\left(\eta+1\right)+log_{2}\left(\eta+2\right)+...+log_{2}\left(N'\right)-\left(N+1\right)log_{2}\left(\eta\right)\\
 & =log_{2}\left(\frac{N'!}{\eta!}\right)-\left(N+1\right)log_{2}\left(\eta\right)\\
 & =N'log_{2}\left(N'\right)-N'log_{2}e+O\left(log_{2}\left(N'\right)\right)-\left(N+1\right)log_{2}\left(\eta\right)-log_{2}\eta!
\end{align*}
 We can now write Eq. \ref{eq:ll-total-instructions} as
\[
total\; instructions\leq\left(N_{x}N_{y}\eta+1\right)
\left[N'log_{2}N' -N'log_{2}e  -\left(N+1\right)log_{2}\eta + N +O\left(log_{2}N'\right)  -log_{2}\eta! +1 \right]
\]
It is clear that the dominant term here is $N'\left(log_{2}\left(N'\right)\right)$and
since $N'$ is simply equal to $N$ offset by a constant, we can conclude
the overall complexity is 
\[
O\left(N_{x}N_{y}Nlog_{2}\left(N\right)\right)
\]
Compared to the full and adaptive step algorithms analyzed in this section, the linked list implementation is clearly
the most efficient in terms of how execution time scales as number of timesteps
$N$, grows large. An added benefit (as noted in \cite{Paper0}, Section 4) is the lowered memory requirements from $O\left(N\right)$ (data stored for all $N$ time steps) to $O\left(log_{2}N\right)$, since we no longer require keeping all time points in the history of the simulation.

As before, we also note that in the one-dimensional analogous case of the linked list pseudocode, the loop iterating over spatial grid points can be vectorized.

\subsection{Simulation Results}
\begin{figure}
\begin{center}
\includegraphics[scale=0.5]{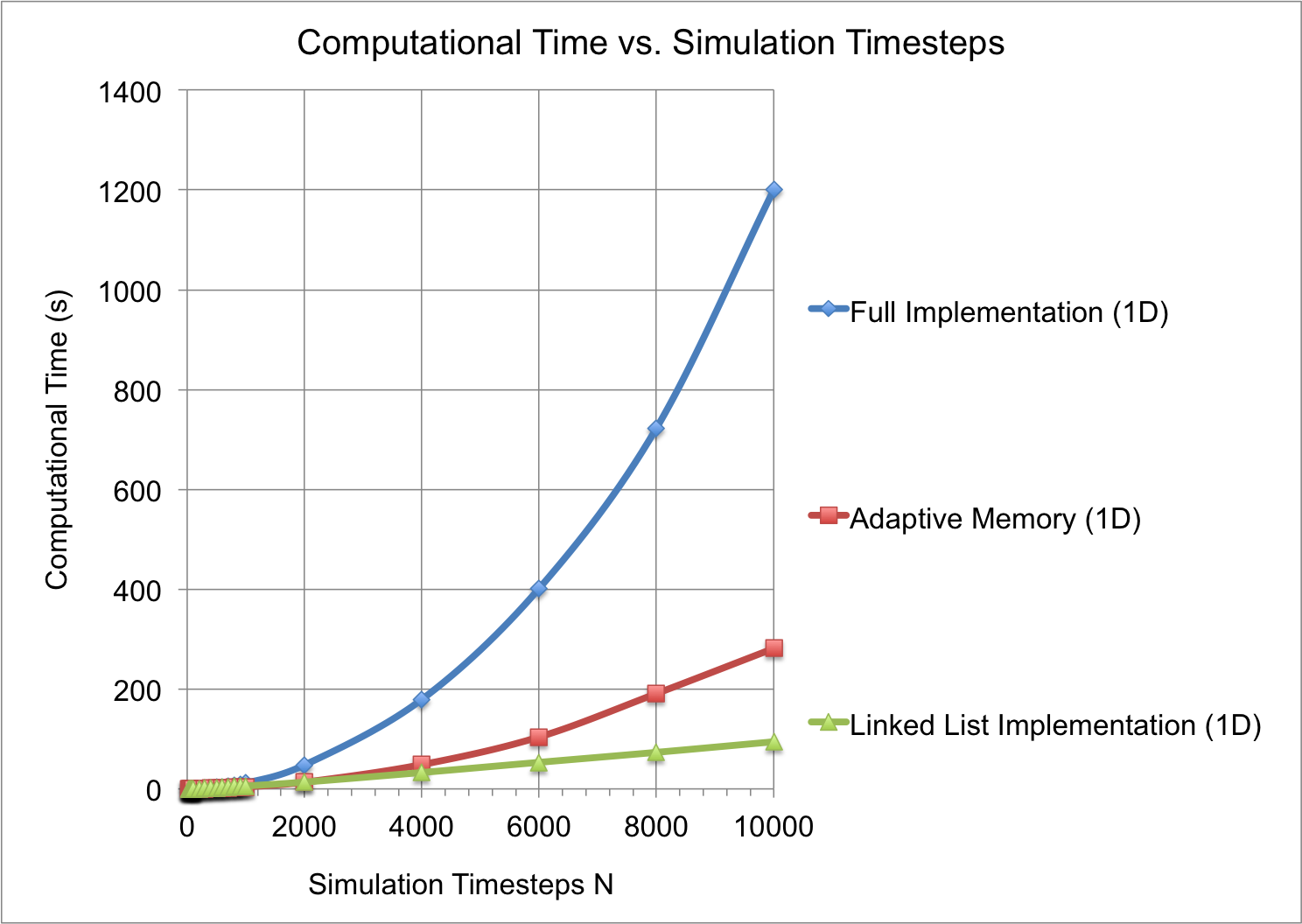}
\caption[Simulated data verifying theoretical complexity results.]
{\footnotesize Simulated data verifying theoretical
complexity results. For all simulations we set $\gamma=0.8,\; D=1 \frac{units^{2}}{s^{\gamma}},\; dt=0.1s,\; dx=0.6875\; units$.
We used adaptive step adaptive step parameter $a$ and $\eta$ values of $20$.}
\label{fig:complexity_simulated}
\end{center}
\end{figure}

We verified our theoretical complexity results with simulated data.
In Fig. \ref{fig:complexity_simulated} we show computation run times
for simulations of various number of steps $N$ and demonstrate the
complexity of the three numerical methods explored in this section (using one-dimensional versions for simplicity). For all algorithms
we set $\gamma=0.8,\; D=1\frac{units^{2}}{s^{\gamma}},\; dt=0.1s,\; dx=0.6875\; units$.
We used adaptive step parameter $a$ and $\eta$ values of
$20$. 

Trendlines were fit to the data using Matlab's curve fitting toolbox.
For the data corresponding to the full and adaptive memory implementations,
the best fit lines were polynomials of second order (modeled with
three coefficients). For the linked list implementation, the best fit trendline with minimum number of coefficients was a $Nlog_{2}N$ function.
\begin{itemize}
\item full implementation simulation time = $p_{1}N^{2}+p_{2}N+p_{3}$ with
$p_{1}=1.106e-5,\; p_{2}=5.5e-4,\; p_{3}=-0.0423$ and goodness of
fit measure $R^{2}=0.9999$
\item adaptive memory simulation time = $p_{1}N^{2}+p_{2}N+p_{3}$ with
$p_{1}=2.671e-5,\; p_{2}=1.441e-3,\; p_{3}=-0.05873$ and goodness
of fit measure $R^{2}=0.9999$
\item linked list simulation time = $p_{1}Nlog_{2}\left(N\right)$ with
$p_{1}=7.047e-4$ and goodness of fit measure $R^{2}=0.9992$.
This is the simplest trendline that was found to be a good fit. There
were other functions that fit the data for the linked list
implementation but required additional coefficients.
\end{itemize}
The data fitting and simulated data in Fig. \ref{fig:complexity_simulated}
thus supports the theoretical analysis done in this section.

It is worth mentioning that there is additional overhead with the linked list implementation, involving the actual handling of the linked list (e.g. deleting and inserting new nodes, etc.). For this reason, for shorter simulations, the linked list algorithm may actually take longer to execute than the full or adaptive memory implementations. It is only when simulation length exceeds a certain point (and this threshold is determined by the details of the simulation, including parameter values) that the $NlogN$ time scaling complexity become more noticeable and the advantages of the linked list approach become more evident.

\section{Error}
\begin{figure}
\begin{center}
\includegraphics[width=6in]{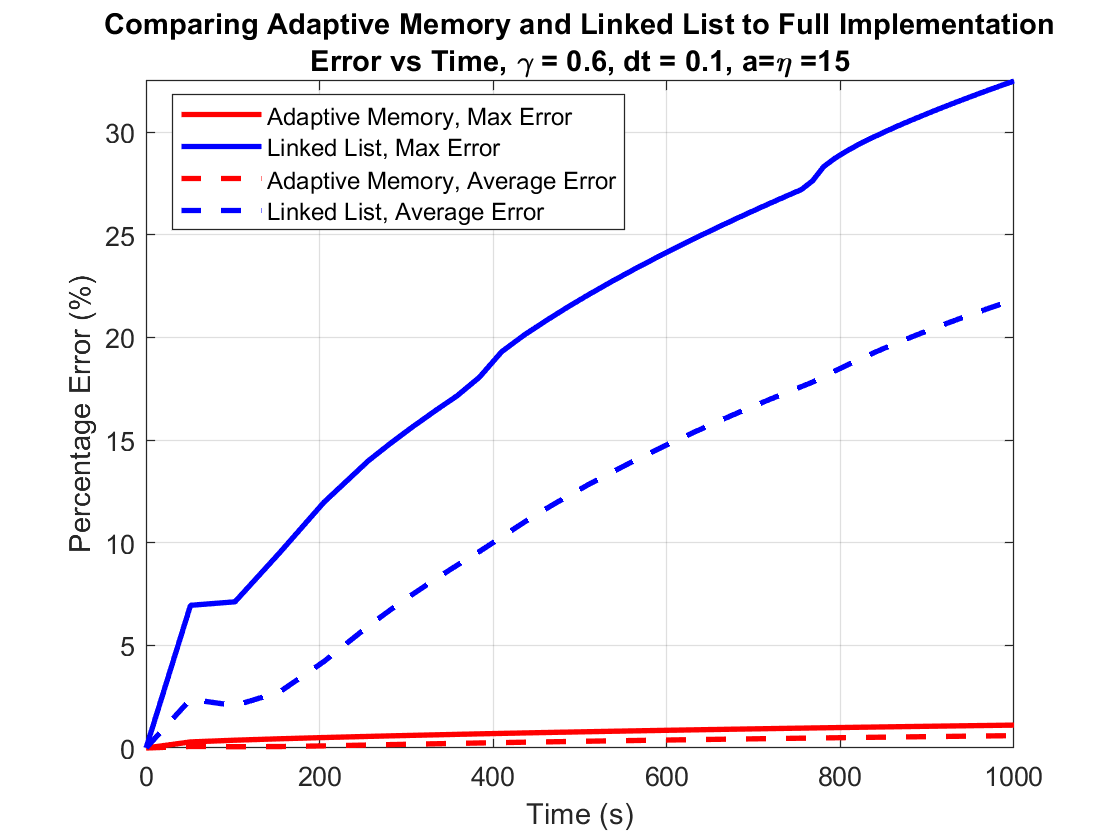}
\caption[Comparing error of adaptive memory and linked list algorithms, using the full implementation as the reference.]
{\footnotesize Comparing error of adaptive memory and linked list algorithms, using the full implementation as the reference. For simulations, the following parameters were used: $\gamma=0.6,\;\alpha=\beta=50\;\frac{units^{2}}{s^{\gamma}},\; dt=0.1s,\;\Delta x=\Delta y=8.64\; units,\; N_{x}=N_{y}=20, \; a = \eta = 15$. We include both maximum error at each timepoint, and average error over the whole solution, at each timepoint.}
\label{fig:AMllvsFull_g0.6}
\end{center}
\end{figure}

In the last section we found that the algorithmic complexity for the
linked list implementation is generally more computationally efficient than the full and adaptive memory implementations, with the exception of shorter simulation lengths where the overhead of dealing with linked list operations, may overwhelm the algorithmic scaling in time. However, with numerical
algorithms, considerable advantages in one area usually come with
some drawback or cost, in another. In this case, the accuracy of the
linked list implementation is the drawback. In Fig. \ref{fig:AMllvsFull_g0.6}
we compare the error of the adaptive memory and linked list approaches, with the full implementation. We include both maximum error at each timepoint, and average error over the whole solution, at each timepoint. The reason we compare with the full implementation is because there is no closed form solution for the 2D fractional diffusion equation except for very specific circumstances, and these two algorithms were derived in part based on the  basic finite difference scheme used in the full implementation. In Fig. \ref{fig:AMllvsFull_g0.6} we see that the adaptive memory algorithm using the given parameters, converges over time, but the error with the linked list approach increases quickly and in an unbounded manner. There is some flexibility on setting parameters, and the larger $\eta$ is, the more accurate the simulation is, but the longer it takes (similar to the adaptive step parameter $a$). We leave as an open problem the analysis of the error resulting from the linked implementation, including whether it converges over time and its exact dependence on $\eta$. Even in
the case that error does not ever converge, one can set the value of $\eta$
depending on the required length of the simulation and how much
error one is willing to tolerate during the simulation. 

In Fig. \ref{fig:AMllvsFull_g0.6_comp_times} we also include a comparison of computational times for the simulations used in Fig. \ref{fig:AMllvsFull_g0.6}. It is clear that the adaptive memory and linked list approaches are much faster than the full implementation. However, the linked list simulation was not considerably faster than the adaptive step approach. If the simulation time was extended further, we would see this gap widen, however, as the $NlogN$ vs $N^2$ scaling effects came into play. Considering both Fig. \ref{fig:AMllvsFull_g0.6} and \ref{fig:AMllvsFull_g0.6_comp_times}, we may come to the conclusion that in many practical cases, the adaptive memory approach may provide the most reasonable trade-off between execution time and accuracy.

\begin{figure}
\begin{center}
\includegraphics[width=6in]{Figures/error/section5_Paper1_time_comp}
\caption[Comparing execution times between all three algorithms.]
{\footnotesize Comparing execution times between all three algorithms. The following parameters were used (same as Fig. \ref{fig:AMllvsFull_g0.6} $\gamma=0.6,\;\alpha=\beta=50\;\frac{units^{2}}{s^{\gamma}},\; dt=0.1s,\;\Delta x=\Delta y=8.64\; units,\; N_{x}=N_{y}=20, \; a = \eta = 15$.}
\label{fig:AMllvsFull_g0.6_comp_times}
\end{center}
\end{figure}

\section{Conclusion}

We have successfully characterized the stability of two finite difference
algorithms used in solving the time-fractional diffusion equation.
We provide our rationale for selecting a von Neumann analysis while
also exploring a few alternative interpretations and approaches to
analyzing stability. Using a von Neumann analysis we have developed
bounding expressions for parameters like time step, spatial discretization,
and grid size, that can be used to appropriately set parameters to
ensure accurate simulations that reflect the true solution of the
fractional diffusion equation. Our simulation results verify that
our bounding expressions resulting from our stability analysis are valid and 
accurate. We also note here that the stability analysis of the linked list implementation should be explored, and we leave this as an open problem for further consideration.

We have also successfully characterized the algorithmic complexity of three finite difference algorithms introduced in \cite{Paper0}. We find that the full and adaptive step algorithms have a Big-$O$ complexity of $O\left(N^{2}\right)$ while the linked list scheme performs much better with $O\left(Nlog_{2}N\right)$ complexity. We compared our analytical results with simulated data and find they are in good agreement.

\FloatBarrier
\bibliographystyle{amsplain}
\bibliography{references_new}
\end{document}